\documentclass[10pt,draft]{amsart}
\usepackage{amscd}
\usepackage{amssymb}

\newtheorem{thm}[equation]{Theorem}
\newtheorem{cor}[equation]{Corollary}
\newtheorem{prop}[equation]{Proposition}
\newtheorem{lem}[equation]{Lemma}
\newtheorem{defn}[equation]{Definition}
\newtheorem{rem}[equation]{Remark}

\numberwithin{equation}{section}

\newcommand{\arr}{\rightarrow}

\newcommand{\lgarr}{\longrightarrow}

\newcommand{\xarr}{\xrightarrow}

\newcommand{\cat}[1]{\operatorname{\mathsf{#1}}}

\newcommand{\opn}{\operatorname}

\newcommand{\rmitem}[1]{\item[\text{\textup{(#1)}}]}

\newcommand{\mcal}[1]{\mathcal{#1}}

\newcommand{\mrm}[1]{\mathrm{#1}}
\newcommand{\mbb}[1]{\mathbb{#1}}

\newcommand{\bsym}[1]{\boldsymbol{#1}}

\newcommand{\Hom}{\operatorname{Hom}}
\newcommand{\End}{\operatorname{End}}
\newcommand{\ten }{\otimes}
\newcommand{\dhom}{\operatorname{Hom}^{\centerdot}}
\newcommand{\dten}{\overset{\centerdot}{\otimes}}
\newcommand{\rhom}{\boldsymbol{R}\operatorname{Hom}^{\centerdot}}
\newcommand{\lten}{\overset{\centerdot}{\otimes}{}^{\boldsymbol{L}}}
\newcommand{\hclim}{\underset{\lgarr}{\cat{hocolim}}}

\title[Recollement and Tilting Complexes]
{Recollement and Tilting Complexes}
\author{Jun-ichi Miyachi}
\date{January 31st, 2002}
\address{J. Miyachi: Department of Mathematics, Tokyo Gakugei
University, Koganei-shi, Tokyo, 184-8501, Japan}
\email{miyachi@u-gakugei.ac.jp}
\subjclass{16G99, 18E30, 18G35}

\begin{document}

\begin{abstract}
First, we study recollement of a derived category of
unbounded complexes of modules induced by a partial tilting complex.
Second, we give equivalent conditions for $P^{\centerdot}$ to be a 
recollement tilting complex, that is,  a tilting complex which induces
an equivalence between recollements $\{\cat{D}_{A/AeA}(A), \\ \cat{D}(A), \cat{D}(eAe)\}$ 
and $\{\cat{D}_{B/BfB}(B), \cat{D}(B), \cat{D}(fBf)\}$, where $e$, $f$ are idempotents of $A$, $B$,
respectively.
In this case, there is an unbounded bimodule complex $\varDelta^{\centerdot}_{T}$ which induces
an equivalence between $\cat{D}_{A/AeA}(A)$ and $\cat{D}_{B/BfB}(B)$.
Third, we apply the above to a symmetric algebra $A$.
We show that a partial tilting complex $P^{\centerdot}$ for $A$ of length 2 extends to a tilting complex, 
and that $P^{\centerdot}$ is a tilting complex if and only if the number of indecomposable 
types of $P^{\centerdot}$ is one of $A$.
Finally, we show that for an idempotent $e$ of $A$, a tilting complex for $eAe$
extends to a recollement tilting complex for $A$,
and that its standard equivalence induces an equivalence between 
$\cat{Mod}A/AeA$ and $\cat{Mod}B/BfB$.
\end{abstract}

\maketitle

\setcounter{section}{-1}
\section{Introduction} \label{s0}

The notion of recollement of triangulated categories was introduced by Beilinson, Bernstein and Deligne
in connection with derived categories of sheaves of topological spaces (\cite{BBD}).
In the representation theory, Cline, Parshall and Scott applied this notion to finite dimensional algebras
over a field, and introduced the notion of quasi-hereditary algebras (\cite{CPS}, \cite{PS}).
In quasi-hereditary algebras, idempotents of algebras play an important role.
In \cite{Rd1}, Rickard introduced the notion of tilting complexes as a generalization of 
tilting modules.
Many constructions of tilting complexes have a relation to idempotents of 
algebras (e.g. \cite{Ok}, \cite{RZ}, \cite{HK1}, \cite{HK2}).
We studied constructions of tilting complexes of term length 2 which has an application to 
symmetric algebras (\cite{HKM}).
In the case of algebras of infinite global dimension, we cannot
treat recollement of derived categories of bounded complexes such as one in the case of quasi-hereditary algebras.
In this paper, we study recollement of derived categories of unbounded complexes of modules 
for projective algebras over a commutative ring $k$, and give the conditions that tilting complexes
induce equivalences between recollements induced by idempotents.
Moreover, we give some constructions of tilting complexes over symmetric algebras.

In Section \ref{s2}, for a projective algebra $A$ over a commutative ring $k$,
we study a recollement $\{\mcal{K}_{P}, \cat{D}(A), \cat{D}(B)\}$
of a derived category $\cat{D}(A)$ of unbounded complexes of right $A$-modules 
induced by a partial tilting complex $P^{\centerdot}$, where $B=\End_{\cat{D}(A)}(P^{\centerdot})$.
We show that there exists the triangle $\xi_{V}$ in $\cat{D}(\cat{Mod}A^{\mrm{e}})$ which
induce adjoint functors of this recollement, and that the triangle $\xi_{V}$ is isomorphic to a triangle
which is constructed by a $P^{\centerdot}$-resolution of $A$ in the sense of Rickard 
(Theorem \ref{parttiltrecoll1}, Proposition \ref{unique01}, Corollary \ref{vardelta2}).
In general, this recollement is out of localizations of triangulated categories which Neeman treated 
in \cite{Ne} (Corollary \ref{cpt1}).
Moreover, we study a recollement $\{\cat{D}_{A/AeA}(A), \cat{D}(A), \cat{D}(eAe)\}$ which
is induced by an idempotent $e$ of $A$ (Proposition \ref{idemprecoll1}, Corollary \ref{idemprecoll2}).
In Section \ref{s3},  we study equivalences between recollements which are induced by idempotents.
We give equivalent conditions for $P^{\centerdot}$ to be a tilting complex inducing an equivalence between
recollements $\{\cat{D}_{A/AeA}(A), \cat{D}(A), \cat{D}(eAe)\}$ and
$\{\cat{D}_{B/BfB}(B), \cat{D}(B), \cat{D}(fBf)\}$ (Theorem \ref{recolltilt1}).
We call this tilting complex a recollement tilting complex related to an idempotent $e$.
There are many symmetric properties between algebras $A$ and $B$ for a 
two-sided recollement tilting complex
${}_{B}T^{\centerdot}_{A}$(Corollaries \ref{recolltilt2}, \ref{recolltilt3}).
Moreover, we have an unbounded bimodule complex $\varDelta^{\centerdot}_{T} \in 
\cat{D}(B^{\circ}\ten A)$ which induces an equivalence between 
$\cat{D}_{A/AeA}(A)$ and $\cat{D}_{B/BfB}(B)$.
The complex $\varDelta^{\centerdot}_{T}$ is a compact object in $\cat{D}_{A/AeA}(A)$,
and satisfies properties such as a tilting complex (Propositions \ref{recolltilt5}, \ref{recolltilt6}, 
\ref{Moritaeqv1}, Corollary \ref{bounded}).
In Section \ref{s4}, we study constructions of tilting complexes for a symmetric algebra $A$
over a field.
First, we construct a family of complexes $\{\Theta_{n}^{\centerdot}(P^{\centerdot}, A)\}_{n \geq 0}$
from a partial tilting complex $P^{\centerdot}$, 
and give equivalent conditions for $\Theta_{n}^{\centerdot}(P^{\centerdot}, A)$ to be a tilting complex
(Definition \ref{consttilt}, Theorem \ref{exttilt2}, Corollary \ref{exttilt3}).
As applications, we show that a partial tilting complex $P^{\centerdot}$
of length 2 extends to a tilting complex, and that $P^{\centerdot}$ is 
a tilting complex if and only if the number of indecomposable types of 
$P^{\centerdot}$ is one of $A$ (Corollaries \ref{exttilt4}, \ref{exttilt5}).  
This is a complex version over symmetric algebras of Bongartz's result on classical tilting modules (\cite{Bo}).
Second, for an idempotent $e$ of $A$, by the above construction a tilting complex for $eAe$
extends to a recollement tilting complex $T^{\centerdot}$ related to $e$ (Theorem \ref{extrecoll}).
This recollement tilting complex induces that $A/AeA$ is isomorphic to $B/BfB$ as a ring, and that
the standard equivalence $\rhom_{A}(T^{\centerdot}, -)$ induces 
an equivalence between $\cat{Mod}A/AeA$ and $\cat{Mod}B/BfB$ (Corollary \ref{Moritaeqv2}).
This construction of tilting complexes contains constructions obtained by several authors.

\section{Basic Tools on Projective Algebras} \label{s1}

In this section, we recall basic tools of derived functors in the case of projective algebras
over a commutative ring $k$.  Throughout this section, we deal only with
projective $k$-algebras, that is, $k$-algebras which are projective as $k$-modules.
For a $k$-algebra $A$, we denote by $\cat{Mod}A$ the category of right $A$-modules, and
denote by $\cat{Proj}A$ (resp., $\cat{proj}A$) the full additive subcategory of 
$\cat{Mod}A$ consisting of projective (resp., finitely generated projective) modules.
For an abelian category $\mcal{A}$ and an additive category $\mcal{B}$, we denote by 
$\cat{D}(\mcal{A})$ (resp., $\cat{D}^{+}(\mcal{A})$, $\cat{D}^{-}(\mcal{A})$, $\cat{D}^{\mrm{b}}(\mcal{A})$)
the derived category of complexes of $\mcal{A}$ (resp.,
complexes of $\mcal{A}$ with bounded below cohomologies, 
complexes of $\mcal{A}$ with bounded above cohomologies, 
complexes of $\mcal{A}$ with bounded cohomologies), 
denote by $\cat{K}(\mcal{B})$ (resp., $\cat{K}^{\mrm{b}}(\mcal{B})$)
the homotopy category of complexes (resp.,
bounded complexes) of $\mcal{B}$ (see \cite {RD} for details).
In the case of $\mcal{A} = \mcal{B} = \cat{Mod}A$, we simply write
$\cat{K}^{*}(A)$ and $\cat{D}^{*}(A)$ for
$\cat{K}^{*}(\cat{Mod}A)$ and $\cat{D}^{*}(\cat{Mod}A)$, respectively.
Given a $k$-algebra $A$ we denote by $A^{\circ}$
the opposite algebra, and by $A^{\mrm{e}}$ the enveloping algebra
$A^{\circ}\otimes_k A$.
We denote by ${Res}_{A}:\cat{Mod}B^{\circ}{\otimes}_{k}A \arr \cat{Mod}A$
the forgetful functor, and use the same symbol 
${Res}_{A}:\cat{D}(B^{\circ}{\otimes}_{k}A) \arr \cat{D}(A)$
for the induced  derived functor.
Throughout this paper, we simply write $\ten$ for $\ten_{k}$.

In the case of projective $k$-algebras $A$, $B$ and $C$, using \cite{CE} Chapter IX Section 2,
we don't need to distinguish the derived functor

\[\begin{aligned}
Res_k\circ(\bsym{R}\opn{Hom}^{\centerdot}_{C}): 
\cat{D}(A^{\circ}\ten C)^{\circ}
\times \cat{D}(B^{\circ}\ten C) \arr
\cat{D}(B^{\circ}\ten A) \arr \cat{D}(k) \\
\text{(resp.,}\ Res_k\circ(\lten_{B} ): 
\cat{D}(A^{\circ}\ten B) 
\times \cat{D}(B^{\circ}\ten C) \arr
\cat{D}(A^{\circ}\ten C) \arr \cat{D}(k))
\end{aligned}\]
with the derived functor 
{\small \[\begin{aligned}
\bsym{R}\opn{Hom}^{\centerdot}_{C}\circ((Res_C)^{\circ}\times Res_C): 
\cat{D}(A^{\circ}\ten C)^{\circ}
\times \cat{D}(B^{\circ}\ten C) \arr
\cat{D}(C)^{\circ} \times \cat{D}(C) \arr \cat{D}(k) \\
\text{(resp.,}\ \lten_{B}\circ(Res_B\times Res_{B^{\circ}}) : 
\cat{D}(A^{\circ}\ten B)
\times \cat{D}(B^{\circ}\ten  C) \arr
\cat{D}(B)\times \cat{D}(B^{\circ}) \arr \cat{D}(k))
\end{aligned}\]}\par\noindent
(see \cite{Rd2}, \cite{BN} and\cite{Sp}  for details).
We freely use this fact in this paper.
Moreover, we have the following statements.

\begin{prop} \label{ddotadj}
Let $k$ be a commutative ring, $A, B, C, D$ projective $k$-algebras.
The following hold.
\begin{enumerate}
\item For ${}_{A}U_{B}^{\centerdot} \in \cat{D}(A^{\circ}\ten B),
{}_{B}V_{C}^{\centerdot} \in \cat{D}(B^{\circ}\ten C),
{}_{C}W_{D}^{\centerdot} \in \cat{D}(C^{\circ}\ten D)$,
we have an isomorphism in $\cat{D}(A^{\circ}\ten D)$:
\[
({}_{A}U^{\centerdot}\lten_{B}V^{\centerdot})
\lten_{C}W_{D}^{\centerdot}
\cong
{}_{A}U^{\centerdot}\lten_{B}(V^{\centerdot}
\lten_{C}W_{D}^{\centerdot}).
\]
\item For ${}_{A}U_{B}^{\centerdot} \in \cat{D}(A^{\circ}\ten B),
{}_{D}V_{C}^{\centerdot} \in \cat{D}(D^{\circ}\ten C),
{}_{A}W_{C}^{\centerdot} \in \cat{D}(D^{\circ}\ten C)$,
we have an isomorphism in $\cat{D}(B^{\circ}\ten D)$:
\[
\bsym{R}\opn{Hom}_{A}^{\centerdot}
({}_{A}U_{B}^{\centerdot},\bsym{R}\opn{Hom}_{C}^{\centerdot}
({}_{D}V_{C}^{\centerdot},{}_{A}W_{C}^{\centerdot})) \cong
\bsym{R}\opn{Hom}_{C}^{\centerdot}
({}_{D}V_{C}^{\centerdot},\bsym{R}\opn{Hom}_{A}^{\centerdot}
({}_{A}U_{B}^{\centerdot},{}_{A}W_{C}^{\centerdot})).
\]
\item For ${}_{A}U_{B}^{\centerdot} \in \cat{D}(A^{\circ}\ten B),
{}_{B}V_{C}^{\centerdot} \in \cat{D}(B^{\circ}\ten C),
{}_{D}W_{C}^{\centerdot} \in \cat{D}(D^{\circ}\ten C)$,
we have an isomorphism in $\cat{D}(D^{\circ}\ten A)$:
\[
\bsym{R}\opn{Hom}_{C}^{\centerdot}
({}_{A}U^{\centerdot}\lten_{B}V_{C}^{\centerdot},
{}_{D}W_{C}^{\centerdot}) \cong
\bsym{R}\opn{Hom}_{B}^{\centerdot}({}_{A}U_{B}^{\centerdot},
\bsym{R}\opn{Hom}_{C}^{\centerdot}({}_{B}V_{C}^{\centerdot},{}_{D}W_{C}^{\centerdot})).
\]
\item For ${}_{A}U_{B}^{\centerdot} \in \cat{D}(A^{\circ}\ten B),
{}_{B}V_{C}^{\centerdot} \in \cat{D}(B^{\circ}\ten C),
{}_{A}W_{C}^{\centerdot} \in \cat{D}(A^{\circ}\ten C)$,
we have an isomorphism in $\cat{D}(k)$:
\[
\bsym{R}\opn{Hom}_{A^{\circ}\ten C}^{\centerdot}
({}_{A}U^{\centerdot}\lten_{B}V_{C}^{\centerdot},{}_{A}W_{C}^{\centerdot}) \cong
\bsym{R}\opn{Hom}_{A^{\circ}\ten B}^{\centerdot}
({}_{A}U_{B}^{\centerdot},\bsym{R}\opn{Hom}_{C}^{\centerdot}
({}_{B}V_{C}^{\centerdot},{}_{A}W_{C}^{\centerdot})).
\]
\item For ${}_{A}U_{B}^{\centerdot} \in \cat{D}(A^{\circ}\ten B),
{}_{B}V_{C}^{\centerdot} \in \cat{D}(B^{\circ}\ten C),
{}_{A}W_{C}^{\centerdot} \in \cat{D}(A^{\circ}\ten C)$,
we have a  commutative diagram:
{\small \[\begin{array}{ccc}
\Hom_{\cat{D}(A^{\circ}\ten C)}
({}_{A}U^{\centerdot}\lten_{B}V_{C}^{\centerdot},{}_{A}W_{C}^{\centerdot})
& \xarr{\sim}
& \Hom_{\cat{D}(A^{\circ}\ten B)}
({}_{A}U_{B}^{\centerdot},\bsym{R}\opn{Hom}_{C}^{\centerdot}
({}_{B}V_{C}^{\centerdot},{}_{A}W_{C}^{\centerdot})) \\
Res_C \downarrow & & \downarrow Res_B \\
\Hom_{\cat{D}(C)}
(U^{\centerdot}\lten_{B}V_{C}^{\centerdot},W_{C}^{\centerdot}) 
& \xarr{\sim}
& \Hom_{\cat{D}(B)}
(U_{B}^{\centerdot},\rhom_{C}
({}_{B}V_{C}^{\centerdot},W_{C}^{\centerdot})) ,
\end{array}\]}\par\noindent
where all horizontal arrows are isomorphisms induced by 3 and 4.
Equivalently, we don't need to distinguish the adjunction arrows induced by
${}_{B}V_{C}^{\centerdot}$ (see \cite{CW}, IV, 7).
\end{enumerate}
\end{prop}

\begin{defn} \label{dfn:bpfcpx}
A complex $X^{\centerdot} \in \cat{D}(A)$ is called
a perfect complex if $X^{\centerdot}$
is isomorphic to a complex of $\cat{K}^{\mrm{b}}(\cat{proj}A)$
in $\cat{D}(A)$.
We denote by $\cat{D}(A)_{\mrm{perf}}$
the triangulated full subcategory of $\cat{D}(A)$
consisting of perfect complexes.
A bimodule complex $X^{\centerdot} \in 
\cat{D}(B^{\circ}{\otimes}_{k}A)$ 
is called a  biperfect complex if $Res_{A}(X^{\centerdot})
\in \cat{D}(A)_{\mrm{perf}}$
and if $Res_{B^{\circ}}(X^{\centerdot})
\in \cat{D}(B^{\circ})_{\mrm{perf}}$.

For an object $C$ of a triangulated category $\mcal{D}$,
$C$ is called a compact object in $\mcal{D}$ if 
$\Hom_{\mcal{D}}(C,-)$ commutes with arbitrary coproducts
on $\mcal{D}$.
\end{defn}

For a complex $X^{\centerdot} = (X^{i}, d^{i})$, we define 
the following truncations:
\[\begin{aligned}
{\sigma}_{\leq n}X^{\centerdot} & : \ldots \arr X^{n-2} \arr X^{n-1}
\arr \opn{Ker}d^{n} \arr 0 \arr \ldots ,\\
{\sigma '}_{\geq n}X^{\centerdot} & : \ldots \arr 0 \arr \opn{Cok}d^{n-1}
\arr X^{n+1}\arr X^{n+2}\arr \ldots .
\end{aligned}\]
The following characterization of perfect complexes is well known
(cf. \cite{Rd1}).  For the convenience of the reader, we give a simple proof.

\begin{prop} \label{pfcp}
For $X^{\centerdot} \in \cat{D}(A)$,
the following are equivalent.
\begin{enumerate}
\item $X^{\centerdot}$ is a perfect complex.
\item $X^{\centerdot}$ is a compact object in $\cat{D}(A)$.
\end{enumerate}
\end{prop}

\begin{proof}
1 $\Rightarrow$ 2.  It is trivial, because we have isomorphisms:
\[\begin{aligned}
\Hom_{\cat{D}(A)}(X^{\centerdot}, - ) & \cong 
R^{0}\dhom_{A}(X^{\centerdot}, -) \\
& \cong \opn{H}^{0}(-\lten_{A}\rhom_{A}(X^{\centerdot},A)).
\end{aligned}\]

2 $\Rightarrow$ 1.
According to \cite{BN} or \cite{Sp}, there is a complex $P^{\centerdot}: \ldots \to 
P^{n-1} \xarr{d^{n-1}} P^{n} \to \ldots \in \cat{K}(\cat{Proj}A)$
such that 
\begin{enumerate}
\rmitem{a} $P^{\centerdot} \cong X^{\centerdot}$ in $\cat{D}(A)$,
\rmitem{b} $\Hom_{\cat{K}(A)}(P^{\centerdot}, - )  \cong
\Hom_{\cat{D}(A)}(P^{\centerdot}, - )$.
\end{enumerate}
Consider the complex $C^{\centerdot}: \ldots \xarr{0} \opn{Cok}d^{n-1} \xarr{0} \ldots$,
then it is easy to see that 
$C^{\centerdot}=$ the coproduct $\coprod_{n \in \mbb{Z}}\opn{Cok}d^{n-1}[-n] = $
the product $\prod_{n \in \mbb{Z}}\opn{Cok}d^{n-1}[-n]$,
that is the biproduct ${\bigoplus}_{n \in \mbb{Z}}\opn{Cok}d^{n-1}[-n]$ of $\opn{Cok}d^{n-1}[-n]$.
Since we have isomorphisms in $\cat{Mod}k$:
\[\begin{aligned}
{\coprod}_{n \in \mbb{Z}}\opn{Hom}_{\cat{K}(A)}
(P^{\centerdot},\opn{Cok}d^{n-1}[-n]) & \cong
\opn{Hom}_{\cat{K}(A)}(P^{\centerdot},
{\bigoplus}_{n \in \mbb{Z}}\opn{Cok}d^{n-1}[-i]) \\
& \cong
{\prod}_{n \in \mbb{Z}}\opn{Hom}_{\cat{K}(A)}
(P^{\centerdot},\opn{Cok}d^{n-1}[-n]) ,
\end{aligned}\]
it is easy to see 
$\opn{Hom}_{\cat{K}(A)}
(P^{\centerdot},\opn{Cok}d^{n-1}[-n]) =0$
for all but finitely many $n \in \mbb{Z}$.
Then there are $m \leq n$ such that $P^{\centerdot} \cong 
{\sigma '}_{\geq m}{\sigma}_{\leq n}P^{\centerdot}$ and
${\sigma '}_{\geq m}{\sigma}_{\leq n}P^{\centerdot} \in 
\cat{K}^{\mrm{b}}(\cat{Proj}A)$.
According to \cite{Rd1} Proposition 6.3, we complete the proof.
\end{proof}

\begin{defn} \label{parttilt}
We call a complex $X^{\centerdot} \in \cat{D}(A)$ a partial tilting complex if 
\begin{enumerate}
\rmitem{a}  $X^{\centerdot} \in \cat{D}(A)_{\mrm{perf}}$,
\rmitem{b}  $\Hom_{\cat{D}(A)}(X^{\centerdot}, X^{\centerdot}[n])=0$ for all $n \not= 0$.
\end{enumerate}
\end{defn}

\begin{defn} \label{assbi}
Let $X^{\centerdot} \in \cat{D}(A)$ be a partial tilting complex,
and $B=\End_{\cat{D}(A)}(X^{\centerdot})$.
According to \cite{Ke} Theorem, there exists a unique bimodule complex
$V^{\centerdot} \in \cat{D}(B^{\circ}\ten A)$ up to isomorphism
such that
\begin{enumerate}
\rmitem{a}  there is an isomorphism 
$\phi:X^{\centerdot} \xarr{\sim} Res_{A}V^{\centerdot}$ in $\cat{D}(A)$
such that $\phi f = \lambda_{B}(f) \phi$ for any $f \in \End_{\cat{D}(A)}(X^{\centerdot})$,
where $\lambda_{B}: B \to \End_{\cat{D}(A)}(V^{\centerdot})$ is the left multiplication morphism.
\end{enumerate}
We call $V^{\centerdot}$ the associated bimodule complex of $X^{\centerdot}$.
In this case, the left multiplication morphism
$\lambda_{B}: B \arr \rhom_{A}(V^{\centerdot},V^{\centerdot})$
is an isomorphism in $\cat{D}(\cat{Mod}B^{\mrm{e}})$.
\end{defn}

Rickard showed that for a tilting complex $P^{\centerdot}$ in $\cat{D}(A)$ 
with $B=\opn{End}_{\cat{D}(A)}(P^{\centerdot})$, there exists a two-sided
tilting complex ${}_{B}T_{A}^{\centerdot} \in \cat{D}(B^{\circ}\ten A)$ (\cite{Rd2}).

\begin{defn} \label{dfn:2tilt2}
A bimodule complex ${}_{B}T_{A}^{\centerdot} \in 
\cat{D}(B^{\circ}{\otimes}_{k}A)$
is called a {\it two-sided tilting complex} provided that
\begin{enumerate}
\rmitem{a} ${}_{B}T_{A}^{\centerdot}$ is a biperfect complex.
\rmitem{b} There exists a biperfect complex ${}_{A}T_{B}^{\vee \centerdot}$
such that \\
(b1) ${}_{B}T^{\centerdot}\lten_{A}
T_{B}^{\vee \centerdot} \cong B$ in 
$\cat{D}(B^{\mrm{e}})$, \\
(b2) ${}_{A}T^{\vee \centerdot}\lten_{B}
T_{A}^{\centerdot} \cong A$ in 
$\cat{D}(A^{\mrm{e}})$. 
\end{enumerate}

We call ${}_{A}T_{B}^{\vee \centerdot}$ the {\it inverse} of 
${}_{B}T_{A}^{\centerdot}$. 
\end{defn}

\begin{prop}[\cite{Rd2}] \label{tilt1}
For a two-sided tilting complex ${}_{B}T_{A}^{\centerdot} \in \cat{D}(B^{\circ}\ten A)$,
the following hold.
\begin{enumerate}
\item We have isomorphisms in $\cat{D}(A^{\circ}\ten B)$:
\[\begin{aligned}
{}_{A}T_{B}^{\vee \centerdot} & \cong 
\bsym{R}\opn{Hom}_{A}^{\centerdot}(T,A) \\
& \cong \bsym{R}\opn{Hom}_{B}^{\centerdot}(T,B).
\end{aligned}\]
\item $\bsym{R}\opn{Hom}_{A}^{\centerdot}(T^{\centerdot},-) \cong 
- \lten_{A}
T^{\vee \centerdot}: \cat{D}^{*}(A) \arr \cat{D}^{*}(B)$ is
a triangle equivalence, and has 
$\bsym{R}\opn{Hom}_{B}^{\centerdot}(T^{\vee \centerdot},-) \cong 
- \lten_{B}
T^{\centerdot}: \cat{D}^{*}(B) \arr \cat{D}^{*}(A)$
as a quasi-inverse, where $*=$ nothing, $+,- , \mrm{b}$.
\end{enumerate}
\end{prop}

In the case of projective $k$-algebras, by \cite{Rd2}
we have also the following result (see also Lemma \ref{pfbi}).

\begin{prop} \label{tilt2}
For a bimodule complex ${}_{B}T_{A}^{\centerdot}$,
the following are equivalent.
\begin{enumerate}
\item ${}_{B}T_{A}^{\centerdot}$ is a two-sided tilting complex.
\item ${}_{B}T_{A}^{\centerdot}$ satisfies that \\
(a) ${}_{B}T_{A}^{\centerdot}$ is a biperfect complex, \\
(b) the right multiplication morphism
$\rho_{A}:A \arr \bsym{R}\opn{Hom}_{B}^{\centerdot}
(T^{\centerdot},T^{\centerdot})$
is an isomorphism in $\cat{D}(\cat{Mod}A^{\mrm{e}})$, \\
(c) the left multiplication morphism
$\lambda_{B}:B \arr \bsym{R}\opn{Hom}_{A}^{\centerdot}
(T^{\centerdot},T^{\centerdot})$
is an isomorphism in $\cat{D}(\cat{Mod}B^{\mrm{e}})$.
\end{enumerate}
\end{prop}

\section{Recollement and Partial Tilting Complexes} \label{s2}

In this section, we study recollements of a derived category $\cat{D}(A)$ induced by
a partial tilting complex $P^{\centerdot}_{A}$ and induced by an idempotent $e$ of $A$.
Throughout this section, all algebras are projective algebras over a commutative ring $k$.

\begin{defn} \label{st-t0}
Let $\mcal{D}, \mcal{D}''$ be triangulated categories, 
and $j^{*}:\mcal{D} \to \mcal{D}''$ a $\partial$-functor.
If $j^{*}$ has a fully faithful right (resp., left) adjoint $j_{*}:\mcal{D}'' \to \mcal{D}$
(resp., $j_{!}:\mcal{D}'' \to \mcal{D}$), then $\{\mcal{D}, \mcal{D}''; j^{*},j_{*}\}$ 
(resp., $\{\mcal{D}, \mcal{D}''; j_{!},j^{*}\}$) is called
a localization (resp., colocalization) of $\mcal{D}$.  Moreover, if $j^{*}$ has a 
fully faithful right adjoint $j_{*}:\mcal{D}'' \to \mcal{D}$ and a fully faithful left adjoint
$j_{!}:\mcal{D}'' \to \mcal{D}$, then $\{\mcal{D}, \mcal{D}''; j_{!}, j^{*},j_{*}\}$ is 
called a bilocalization of $\mcal{D}$.

For full subcategories $\mcal{U}$ and $\mcal{V}$ of $\mcal{D}$,
$(\mcal{U}, \mcal{V})$ is called a {\it stable $t$-structure} in $\mcal{D}$ provided that
\begin{enumerate}
\rmitem{a}  $\mcal{U}$ and $\mcal{V}$ are stable for translations.  
\rmitem{b}  $\opn{Hom}_{\mcal{D}}(\mcal{U}, \mcal{V}) = 0$.
\rmitem{c}  For every $X  \in \mcal{D}$, there exists a triangle $U \arr X \arr V \arr U[1]$
with $U \in \mcal{U}$ and $V \in \mcal{V}$.
\end{enumerate}
\end{defn}

We have the following properties.

\begin{prop}[\cite{BBD}, cf. \cite{Mi}] \label{t-st1}
Let $(\mcal{U}, \mcal{V})$ be a stable $t$-structure in a triangulated category $\mcal{D}$,
and let $U \arr X \arr V \arr U[1]$ and $U' \arr X' \arr V' \arr U'[1]$ be triangles
in $\mcal{D}$ with $U, U' \in \mcal{U}$ and $V, V' \in \mcal{V}$.
For any morphism $f : X \to X'$, there exist a unique $f_{\mcal{U}} : U \to U'$ and
a unique $f_{\mcal{V}} : V \to V'$ which induce a morphism of triangles:
\[\begin{CD}
U @>>> X @>>> V @>>> U[1] \\
@Vf_{\mcal{U}}VV @VVf V @VVf_{\mcal{V}}V @VVf_{\mcal{U}}[1]V \\
U' @>>> X' @>>> V' @>>> U'[1] .
\end{CD}\]
In particular, for any $X \in \mcal{D}$, the above $U$ and $V$ are uniquely determined up to
isomorphism.
\end{prop}

\begin{prop}[\cite{Mi}] \label{bilocal}
The following hold.
\begin{enumerate}
\item If $\{\mcal{D}, \mcal{D}''; j^{*},j_{*}\}$ 
(resp., $\{\mcal{D}, \mcal{D}''; j_{!},j^{*}\}$) is a localization (resp., a colocalization) of $\mcal{D}$, then
$(\opn{Ker}j^{*}, \opn{Im}j_{*})$ (resp., $(\opn{Im}j_{!}, \opn{Ker}j^{*})$) is a stable $t$-structure.
In this case, the adjunction arrow $\bsym{1}_{\mcal{D}} \to j_{*}j^{*}$
(resp., $j_{!}j^{*} \to \bsym{1}_{\mcal{D}}$) implies triangles
\[\begin{aligned}
U \to X \to j_{*}j^{*}X \to U[1] \\
\text{ (resp.,}\ j_{!}j^{*}X \to X \to V \to X[1] ) 
\end{aligned}\]
with $U \in \opn{Ker}j^{*}$, $j_{*}j^{*}X \in \opn{Im}j_{*}$
 (resp., $j_{!}j^{*}X \in \opn{Im}j_{!}$, $V \in \opn{Ker}j^{*}$) for all $X \in \mcal{D}$.
\item  If $\{\mcal{D}, \mcal{D}'';j_{!},j^{*},j_{*}\}$ is a bilocalization of $\mcal{D}$, then
the canonical embedding $i_{*}:\opn{Ker}j^{*} \to \mcal{D}$ has  a right adjoint 
$i^{!}:\mcal{D} \to \opn{Ker}j^{*}$and a left adjoint $i^{*}:\mcal{D} \to \opn{Ker}j^{*}$ such that 
$\{\opn{Ker}j^{*}, \mcal{D}, \mcal{D}''; i^{*}, i_{*},i^{!}, j_{!},j^{*},j_{*}\}$ is a recollement
in the sense of \cite{BBD}.
\item If $\{\mcal{D}', \mcal{D}, \mcal{D}''; i^{*}, i_{*},i^{!}, j_{!},j^{*},j_{*}\}$ is a recollement,
then  $\{\mcal{D}, \mcal{D}'';j_{!},j^{*},j_{*}\}$ is a bilocalization of
$\mcal{D}$.
\end{enumerate}
\end{prop}

\begin{prop}[\cite{BBD}] \label{t-st2}
Let $\{\mcal{D}', \mcal{D}, \mcal{D}''; i^{*}, i_{*},i^{!}, j_{!},j^{*},j_{*}\}$ be a recollement,
then \\ $(\opn{Im}i_{*}, \opn{Im}j_{*})$ and $(\opn{Im}j_{!}, \opn{Im}i_{*})$ are stable $t$-structures
in $\mcal{D}$. Moreover, the adjunction arrows 
$\alpha:i_{*}i^{!} \to \bsym{1}_{\mcal{D}}$, $\beta:\bsym{1}_{\mcal{D}} \to j_{*}j^{*}$,
$\gamma:j_{!}j^{*} \to \bsym{1}_{\mcal{D}}$, $\delta:\bsym{1}_{\mcal{D}} \to i_{*}i^{*}$
imply triangles in $\mcal{D}$:
\[\begin{aligned}
i_{*}i^{!}X \xarr{\alpha_{X}} X \xarr{\beta_{X}} j_{*}j^{*}X \to i_{*}i^{!}X[1] , \\
j_{!}j^{*}X \xarr{\gamma_{X}} X \xarr{\delta_{X}} i_{*}i^{*}X \to j_{!}j^{*}X[1],
\end{aligned}\]
for any $X \in \mcal{D}$.
\end{prop}

By Definition \ref{st-t0}, we have the following properties.

\begin{cor} \label{t-st3}
Under the condition of Proposition \ref{t-st2}, the following hold for $X \in \mcal{D}$.
\begin{enumerate}
\item  $i_{*}i^{!}X \cong X$ (resp., $X \cong j_{*}j^{*}X$) in $\mcal{D}$
if and only if $\alpha_{X}$ (resp., $\beta_{X}$) is an isomorphism.
\item  $j_{!}j^{*}X \cong X$ (resp., $X \cong i_{*}i^{*}X$) in $\mcal{D}$
if and only if $\gamma_{X}$ (resp., $\delta_{X}$) is an isomorphism.
\end{enumerate}
\end{cor}

For $X \in \cat{Mod}C^{\circ}\ten A$, $Q \in \cat{Mod}B^{\circ}\ten A$,
let 
\[
\tau_{Q}(X): X\ten_{A}\Hom_{A}(Q, A) \to \Hom_{A}(Q, X)
\]
be the morphism in $\cat{Mod}C^{\circ}\ten B$ defined by 
$(x{\otimes}f \mapsto (q \mapsto xf(q)))$ for $x \in X, q \in Q$, $f \in \Hom_{A}(Q, A)$.
We have the following functorial isomorphism of derived functors.

\begin{lem} \label{pfbi}
Let $k$ be a commutative ring, $A$, $B$, $C$ projective $k$-algebras,
${}_{B}V^{\centerdot}_{A} \in \cat{D}(B^{\circ}\ten A)$ with
$Res_{A}V^{\centerdot} \in \cat{D}(A)_{\mrm{perf}}$, and
$V^{\star \centerdot}=\rhom_{A}(V^{\centerdot},A)\in \cat{D}(A^{\circ}\ten B)$.
Then we have the ($\partial$-functorial) isomorphism:
\[
\tau_{V}: - \lten_{A}V^{\star \centerdot} \xarr{\sim} \rhom_{A}(V^{\centerdot}, -) 
\]
as derived functors $\cat{D}(C^{\circ}\ten A) \to \cat{D}(C^{\circ}\ten B)$.
\end{lem}

\begin{proof}
It is easy to see that we have a $\partial$-functorial
morphism of derived functors \\
$\cat{D}(C^{\circ}\ten A) \arr \cat{D}(C^{\circ}\ten B)$:
\[
\tau_{V}: - \lten_{A}V^{\star \centerdot} \to \rhom_{A}(V^{\centerdot}, -)  .
\]
Let $P^{\centerdot} \in 
\cat{K}^{\mrm{b}}(\cat{proj}A)$
which has a quasi-isomorphism $P^{\centerdot} \arr 
Res_{A}V^{\centerdot}$.
Then we have a $\partial$-functorial isomorphism
of $\partial$-functors
$\cat{D}(C^{\circ}\ten A) \arr \cat{D}(C^{\circ})$
\[
\tau_{P}:
-\dten_{A}
\opn{Hom}_{A}^{\centerdot}(P^{\centerdot},A)
\xarr{\sim}
\Hom_{A}^{\centerdot}(P^{\centerdot},-) .
\]
Since $Res_{C^{\circ}}\circ\tau_{V} \cong \tau_{P}$
and $\opn{H}^{\centerdot}(\tau_{P})$ is an isomorphism,
$\tau_{V}$ is a $\partial$-functorial isomorphism.
\end{proof}

Concerning adjoints of the derived functor $-\lten_{A}V^{\star \centerdot}$, by direct calculation
we have the following properties.

\begin{lem} \label{pfadj}
Let $k$ be a commutative ring, $A$, $B$, $C$ projective $k$-algebras,
${}_{B}V^{\centerdot}_{A} \in \cat{D}(B^{\circ}\ten A)$ with
$Res_{A}V^{\centerdot} \in \cat{D}(A)_{\mrm{perf}}$, and
${}_{A}V_{B}^{\star \centerdot}=\rhom_{A}(V^{\centerdot},A) \in \cat{D}(A^{\circ}\ten B)$.
Then the following hold.
\begin{enumerate}
\item $\tau_{V}$ induces the adjoint isomorphism:
\[\begin{aligned}
\Phi: 
\Hom_{\cat{D}(C^{\circ}\ten B)}(-, ?\lten_{A}V^{\star \centerdot})
& \xarr{\sim}
\Hom_{\cat{D}(C^{\circ}\ten A)}(-\lten_{B}V^{\centerdot}, ?) .
\end{aligned}\]
Therefore, we get the morphism $\varepsilon_{V}: 
V^{\star \centerdot}\lten_{B}V^{\centerdot} \to A$ in
$\cat{D}(A^{\mrm{e}})$
(resp., $\vartheta_{V}: B \to V^{\centerdot}\lten_{A}V^{\star \centerdot}$ in
$\cat{D}(B^{\mrm{e}})$)
from the adjunction arrow
of $A \in \cat{D}(A^{\mrm{e}})$ (resp., $B \in \cat{D}(B^{\mrm{e}})$).
\item In the adjoint isomorphism of 1, the adjunction arrow 
$-\lten_{A}V^{\star \centerdot}\lten_{B}V^{\centerdot} \to
\bsym{1}_{\cat{D}(C^{\circ}\ten A)}$ 
(resp., $\bsym{1}_{\cat{D}(C^{\circ}\ten B)} \to
-\lten_{B}V^{\centerdot}\lten_{A}V^{\star \centerdot}$)
is isomorphic to
$-\lten_{A}\varepsilon_{V}$ (resp., $-\lten_{B}\vartheta_{V}$).
\item In the adjoint isomorphism: 
\[
\Hom_{\cat{D}(C^{\circ}\ten A)}(-, \rhom_{B}(V^{\star \centerdot},?)) 
\xarr{\sim}
\Hom_{\cat{D}(C^{\circ}\ten B)}(-\lten_{A}V^{\star \centerdot}, ?) ,
\]
the adjunction arrow $\bsym{1}_{\cat{D}(C^{\circ}\ten A)} \to
\rhom_{B}(V^{\star \centerdot}, -\lten_{A}V^{\star \centerdot})$
(resp., \\ $\rhom_{B}(V^{\star \centerdot}, -)\lten_{A}V^{\star \centerdot} \to
\bsym{1}_{\cat{D}(C^{\circ}\ten B)}$)
is isomorphic to
$\rhom_{A}(\varepsilon_{V}, -)$ (resp., $\rhom_{B}(\vartheta_{V}, -)$).
\end{enumerate}
\end{lem}

Let $A$, $B$ be projective algebras over a commutative ring $k$.
For a partial tilting complex $P^{\centerdot} \in \cat{D}(A)$ with 
$B\cong\End_{\cat{D}(A)}(P^{\centerdot})$, let ${}_{B}V_{A}^{\centerdot}$ be the associated bimodule
complex of $P^{\centerdot}$.
By Lemma \ref{pfbi}, we can take
 \[\begin{aligned}
j_{V !} & = - \lten_{B}V^{\centerdot} : \cat{D}(B) \to \cat{D}(A) , \\
j_{V}^{*} & = -\lten_{A}V^{\star \centerdot} \cong \rhom_{A}(V^{\centerdot},-) : \cat{D}(A) \to \cat{D}(B) , \\
j_{V *} & = \rhom_{B}(V^{\star \centerdot},-) : \cat{D}(B) \to \cat{D}(A).
\end{aligned}\]
By Lemma \ref{pfadj}, we get
the triangle $\xi_{V}$ in $\cat{D}(A^{\mrm{e}})$:
\[
V^{\star \centerdot}\lten_{B}V^{\centerdot} \xarr{\varepsilon_{V}} A \xarr{\eta_{V}} 
\varDelta^{\centerdot}_{A}(V^{\centerdot}) \to V^{\star \centerdot}\lten_{B}V^{\centerdot}[1].
\]
Let $\mcal{K}_{P}$ be the full subcategory
of $\cat{D}(A)$ consisting of complexes $X^{\centerdot}$ such that
$\Hom_{\cat{D}(A)}(P^{\centerdot}, X^{\centerdot}[i])=0$ for all $i \in \mbb{Z}$.

\begin{thm} \label{parttiltrecoll1}
Let $A$, $B$ be projective algebras over a commutative ring $k$,
$P^{\centerdot} \in \cat{D}(A)$ a partial tilting complex with 
$B\cong\End_{\cat{D}(A)}(P^{\centerdot})$, and let ${}_{B}V_{A}^{\centerdot}$ be the associated bimodule
complex of $P^{\centerdot}$.
Take
{\small \[\begin{aligned}
j_{V !} & = - \lten_{B}V^{\centerdot} : \cat{D}(B) \to \cat{D}(A) ,
& j_{V}^{*} & = -\lten_{A}V^{\star \centerdot} : \cat{D}(A) \to \cat{D}(B) , \\
j_{V *} & = \rhom_{B}(V^{\star \centerdot}, -) : \cat{D}(B) \to \cat{D}(A) , 
& \ i^{*}_{V} & = - \lten_{A}\varDelta^{\centerdot}_{A}(V^{\centerdot}) : \cat{D}(A) \to 
\mcal{K}_{P}, \\
i_{V *} & =\ \text{the embedding }: \mcal{K}_{P} \to \cat{D}(A) ,
& i^{!}_{V} & = \rhom_{A}(\varDelta^{\centerdot}_{A}(V^{\centerdot}), - ) : 
\cat{D}(A) \to \mcal{K}_{P} ,
\end{aligned}\]}\par\noindent
then $\{\mcal{K}_{P}, \cat{D}(A), \cat{D}(B);  
i^{*}_{V}, i_{V *}, i^{!}_{V}, j_{V !}, j_{V}^{*}, j_{V *}\}$:
\[\begin{array}{ccccc}
\mcal{K}_{P} & \overset{\leftarrow}{\rightleftarrows} &
\cat{D}(A) & \overset{\leftarrow}{\rightleftarrows} &
\cat{D}(B)
\end{array}\]
is a recollement.
\end{thm}

\begin{proof}
Since it is easy to see that $\tau_{V}(V^{\centerdot})\circ \vartheta_{V}$
is the left multiplication morphism $B \to \rhom_{A}(V^{\centerdot}, V^{\centerdot})$,
by the remark of Definition \ref{assbi}, $\vartheta_{V}: B \to 
V^{\centerdot}\lten_{A}V^{\star \centerdot}$ is an isomorphism in $\cat{D}(B^{\mrm{e}})$.
By Lemma \ref{pfadj},
$\{\cat{D}(A), \cat{D}(B); j_{V !}, j_{V}^{*}, j_{V *}\}$ is a bilocalization.
By Proposition \ref{bilocal}, there exist $i^{*}_{V}:\cat{D}(A) \to \mcal{K}_{P}$,
$i_{V *}=$ the embedding $:\mcal{K}_{P} \to \cat{D}(A)$, 
$i^{!}_{V}: \cat{D}(A) \to \mcal{K}_{P}$ such that 
 $\{\mcal{K}_{P}, \cat{D}(A), \cat{D}(B);  
i^{*}_{V}, i_{V *}, i^{!}_{V}, j_{V !}, j_{V}^{*}, j_{V *}\}$ is a recollement.
For $X^{\centerdot} \in \cat{D}(A)$, by Lemma \ref{pfadj}, 
$X^{\centerdot}\lten_{A}\varepsilon_{V}$ is isomorphic to the adjunction arrow
$j_{V !}j_{V}^{*}(X^{\centerdot}) \to X^{\centerdot}$.
Then $X^{\centerdot}\lten_{A}\eta_{V}$ is isomorphic to the adjunction arrow
$X^{\centerdot} \to i_{V *}i_{V}^{*}(X^{\centerdot})$, and hence
we can take $i_{V}^{*}= - \lten_{A}\varDelta^{\centerdot}_{A}(V^{\centerdot})$
by Propositions \ref{t-st1}, \ref{t-st2}.
Similarly, 
we can take $i_{V}^{!}= \rhom_{A}(\varDelta^{\centerdot}_{A}(V^{\centerdot}), -)$.
\end{proof}

In general, the above $\varDelta_{A}^{\centerdot}(V^{\centerdot})$ and 
$\varDelta_{A}^{\centerdot}(e)$ in Proposition \ref{idemprecoll1} are unbounded complexes.
Then, by the following corollary we have unbounded complexes
which are compact objects in $\mcal{K}_{P}$ and in $\cat{D}_{A/AeA}(A)$.  
This shows that recollements of Theorem \ref{parttiltrecoll1} and  Proposition \ref{idemprecoll1}
are out of localizations of triangulated categories which Neeman treated in \cite{Ne}.

\begin{cor} \label{cpt1}
Under the condition Theorem \ref{parttiltrecoll1}, the following hold.
\begin{enumerate}
\item $\mcal{K}_{P}$ is closed under coproducts in $\cat{D}(A)$.
\item For any $X^{\centerdot} \in \cat{D}(A)_{\mrm{perf}}$,
$X^{\centerdot}\lten_{A}\varDelta^{\centerdot}_{A}(V^{\centerdot})$ is a compact object in $\mcal{K}_{P}$.
\end{enumerate}
\end{cor}

\begin{proof}
1. 
Since $P^{\centerdot}$ is a compact object in $\cat{D}(A)$, it is trivial.

2.
Since we have an isomorphism:
\[
\Hom_{\cat{D}(A)}(i^{*}_{V}X^{\centerdot}, Y^{\centerdot})
\cong \Hom_{\cat{D}(A)}(X^{\centerdot}, Y^{\centerdot}) 
\]
for any $Y^{\centerdot} \in \mcal{K}_{P}$, we have the statement.
\end{proof}

\begin{cor} \label{vardelta1}
Let $A$, $B$ be projective algebras over a commutative ring $k$,
$P^{\centerdot} \in \cat{D}(A)$ a partial tilting complex with 
$B\cong\End_{\cat{D}(A)}(P^{\centerdot})$, and let ${}_{B}V_{A}^{\centerdot}$ be the associated bimodule
complex of $P^{\centerdot}$.  Then the following hold.
\begin{enumerate}
\item $\varDelta_{A}^{\centerdot}(V^{\centerdot})
\cong \varDelta_{A}^{\centerdot}(V^{\centerdot})\lten_{A}\varDelta_{A}^{\centerdot}(V^{\centerdot})$
in $\cat{D}(A^{\mrm{e}})$.
\item $\rhom_{A}(\varDelta_{A}^{\centerdot}(V^{\centerdot}),
\varDelta_{A}^{\centerdot}(V^{\centerdot})) \cong \varDelta_{A}^{\centerdot}(V^{\centerdot})$
in $\cat{D}(A^{\mrm{e}})$.
\end{enumerate}
\end{cor}

\begin{proof}
Since $\varDelta_{A}^{\centerdot}(V^{\centerdot})\lten_{A}V^{\star \centerdot}[n] \cong
j^{*}_{V}i_{V *}i^{*}_{V}(A[n]) =O$ for all $n$, $\varDelta_{A}^{\centerdot}(V^{\centerdot})\lten_{A}\eta_{V}$ is 
an isomorphism in $\cat{D}(A^{\mrm{e}})$.  Similarly, since
\[\begin{aligned}
\rhom_{A}(V^{\star \centerdot}\lten_{B}V^{\centerdot}, \varDelta_{A}^{\centerdot}(V^{\centerdot}))[n]
& \cong
\rhom_{B}(V^{\star \centerdot}, \varDelta_{A}^{\centerdot}(V^{\centerdot})\lten_{A}V^{\star \centerdot})[n] \\
& =O
\end{aligned}\]
for all $n$,
$\rhom_{A}(\eta_{V}, \varDelta_{A}^{\centerdot}(V^{\centerdot}))$ is an isomorphism in $\cat{D}(A^{\mrm{e}})$.
\end{proof}

\begin{lem} \label{hlim1}
Let $\mcal{D}$ be a triangulated category.  Then the following hold.
\begin{enumerate}
\item  For morphisms of triangles in $\mcal{D}$ ($n\geq 1$):
\[\begin{CD}
X_{n} @>>> Y_{n} @>>> Z_{n} @>>> X_{n}[1] \\
@VVV @VVV @VVV @VVV \\
X_{n+1} @>>> Y_{n+1} @>>> Z_{n+1} @>>> X_{n+1}[1],
\end{CD}\]
there exists a triangle $\coprod X_{n} \to \coprod X_{n} \to X \to \coprod X_{n}[1]$ such that
we have the following triangle in $\mcal{D}$:
\[
X \to \hclim Y_{n} \to \hclim Z_{n} \to X[1] .
\]
\item For a family of triangles in $\mcal{D}$:
$
C_{n} \to X_{n-1} \to X_{n} \to C_{n}[1] \ (n \geq 1),
$
with $X_{0}=X$, there exists a family of triangles in $\mcal{D}$:
\[
C_{n}[-1] \to Y_{n-1} \to Y_{n} \to C_{n} \ (n \geq 1),
\]
with $Y_{0}=O$, such that we have the following triangle in $\mcal{D}$:
\[
Y \to X \to \hclim X_{n} \to Y[1] ,
\]
where $\coprod Y_{n} \to \coprod Y_{n} \to Y \to \coprod Y_{n}[1]$
is a triangle in $\mcal{D}$.
\end{enumerate}
\end{lem}

\begin{proof}
1.  By the assumption, we have a commutative diagram:
\[\begin{CD}
\coprod X_{n} @>>> \coprod Y_{n} @>>> \coprod Z_{n} @>>> \coprod X_{n}[1] \\
@. @VV1-\text{shift}V @VV1-\text{shift}V @. \\
\coprod X_{n} @>>> \coprod Y_{n} @>>> \coprod Z_{n} @>>> \coprod X_{n}[1] \\
\end{CD}\]
According to \cite{BBD} 9 lemma, we have the statement.

2.
By the octahedral axiom, we have a commutative diagram:
\[\begin{CD}
@. @. C_{n} @= C_{n} \\
@. @. @VVV @VVV \\
Y_{n-1} @>>> X @>>> X_{n-1} @>>> Y_{n-1}[1] \\
@VVV @| @VVV @VVV \\
Y_{n} @>>> X @>>> X_{n} @>>> Y_{n}[1] \\
@. @. @VVV @VVV \\
@. @. C_{n}[1] @= C_{n}[1] ,
\end{CD}\]
where all lines are triangles in $\mcal{D}$.
By 1, we have the statement.
\end{proof}

For an object $M$ in an additive category $\mcal{B}$, 
we denote by $\cat{Add}M$ (resp., $\cat{add}M$) the full subcategory of $\mcal{B}$
consisting of objects which are isomorphic to summands of coproducts 
(resp., finite coproducts) of copies of $M$.

\begin{defn} \label{constadj}
Let $A$ be a projective algebra over a commutative ring $k$,
and $P^{\centerdot} \in \cat{D}(A)$ a partial tilting complex.
For $X^{\centerdot} \in \cat{D}^{-}(A)$, there exists an integer $r$
such that $\Hom_{\cat{D}(A)}(P^{\centerdot}, X^{\centerdot}[r+i])=0$ for all $i >0$.
Let $X_{0}^{\centerdot} = X^{\centerdot}$.  For $n \geq 1$, by induction we construct 
a triangle:
\[
P^{\centerdot}_{n}[n-r-1] \xarr{g_n} X_{n-1}^{\centerdot}
\xarr{h_n} X_{n}^{\centerdot} \arr P^{\centerdot}_{n}[n-r]
\]
as follows.  If $\opn{Hom}_{\cat{D}(A)}(P^{\centerdot}, 
X_{n-1}^{\centerdot}[r-n+1]) = 0$,  then we set $P^{\centerdot}_{n} = O$.  
Otherwise, we take $P^{\centerdot}_{n} \in \cat{Add}P^{\centerdot}$ and a morphism 
$g'_n : P^{\centerdot}_{n}  \arr X_{n-1}^{\centerdot}[r-n+1]$ such 
that $\opn{Hom}_{\cat{D}(A)}(P^{\centerdot},  g'_n)$ is an epimorphism, and let 
$g_n = g'_n[n-r-1]$.
By Lemma \ref{hlim1}, we have triangles:
\[
P^{\centerdot}_{n}[n-r-2] \to Y_{n-1}^{\centerdot}
\to Y_{n}^{\centerdot} \arr P^{\centerdot}_{n}[n-r-1]
\]
and $Y_{0}^{\centerdot}=O$.
Then we define 
$\nabla^{\centerdot}_{\infty}(P^{\centerdot},X^{\centerdot})$ and
$\varDelta^{\centerdot}_{\infty}(P^{\centerdot},X^{\centerdot})$ by
the complex of  Lemma \ref{hlim1} (2) and
$\hclim X_{n}^{\centerdot}$, respectively.
Moreover, we have a triangle:
\[
\nabla^{\centerdot}_{\infty}(P^{\centerdot},X^{\centerdot}) \to
X^{\centerdot} \to
\varDelta^{\centerdot}_{\infty}(P^{\centerdot},X^{\centerdot}) \to
\nabla^{\centerdot}_{\infty}(P^{\centerdot},X^{\centerdot})[1] .
\]
\end{defn}

\begin{lem} \label{unique00}
Let $A$, $B$ be projective algebras over a commutative ring $k$, 
$P^{\centerdot} \in \cat{D}(A)$ a partial tilting complex with $B \cong \End_{\cat{D}(A)}(P^{\centerdot})$,
and ${}_{B}V^{\centerdot}_{A}$ the associated bimodule complex of $P^{\centerdot}$.
For $X^{\centerdot} \in \cat{D}^{-}(A)$, we have an isomorphism of triangles in $\cat{D}(A)$:
\[\begin{CD}
j_{V !}j^{*}_{V}X^{\centerdot} @>>> X^{\centerdot} @>>> i_{V *}i^{*}_{V}X^{\centerdot}
@>>>j_{V !}j^{*}_{V}X^{\centerdot}[1] \\
@VV\wr V @| @VV\wr V @VV\wr V \\
\nabla^{\centerdot}_{\infty}(P^{\centerdot},X^{\centerdot}) @>>>
X^{\centerdot} @>>>
\varDelta^{\centerdot}_{\infty}(P^{\centerdot},X^{\centerdot}) @>>>
\nabla^{\centerdot}_{\infty}(P^{\centerdot},X^{\centerdot})[1] .
\end{CD}\]
\end{lem}

\begin{proof}
By the construction, we have $\Hom_{\cat{D}(A)}(P^{\centerdot}, 
\varDelta^{\centerdot}_{\infty}(P^{\centerdot},X^{\centerdot})[i])=0$
for all $i$, and then 
$\varDelta^{\centerdot}_{\infty}(P^{\centerdot},X^{\centerdot}) \in \opn{Im}i_{V *}$
(see Lemma \ref{exttilt1}).
Since $j_{V !}$ is fully faithful and $P^{\centerdot} \in \opn{Im}j_{V !}$, 
it is easy to see $Y^{\centerdot}_{n} \in \opn{Im}j_{V !}$.  Then
$\nabla^{\centerdot}_{\infty}(P^{\centerdot},X^{\centerdot}) \in \opn{Im}j_{V !}$,
because $j_{V !}$ commutes with coproducts.
By Proposition \ref{t-st1}, we complete the proof.
\end{proof}

\begin{defn} \label{constadj2}
Let $A$ be a projective algebra over a commutative ring $k$,
and $P^{\centerdot} \in \cat{D}(A)$ a partial tilting complex.
Given $X^{\centerdot} \in \cat{D}(A)$, for $n \geq 0$,
we have a triangle:
\[
\nabla_{\infty}^{\centerdot}(P^{\centerdot}, \sigma_{\leq n}X^{\centerdot})
\to
\sigma_{\leq n}X^{\centerdot} \to
\varDelta_{\infty}^{\centerdot}(P^{\centerdot}, \sigma_{\leq n}X^{\centerdot})
\to
\nabla_{\infty}^{\centerdot}(P^{\centerdot}, \sigma_{\leq n}X^{\centerdot})[1] .
\]
According to Lemma \ref{unique00} and Proposition \ref{t-st1}, for $n \geq 0$
we have a morphism of triangles:
{\scriptsize \[\begin{array}{ccccccc}
\nabla_{\infty}^{\centerdot}(P^{\centerdot}, \sigma_{\leq n}X^{\centerdot})
& \to &
\sigma_{\leq n}X^{\centerdot} & \to &
\varDelta_{\infty}^{\centerdot}(P^{\centerdot}, \sigma_{\leq n}X^{\centerdot})
& \to &
\nabla_{\infty}^{\centerdot}(P^{\centerdot}, \sigma_{\leq n}X^{\centerdot})[1] \\
\downarrow & & \downarrow & & \downarrow & & \downarrow  \\
\nabla_{\infty}^{\centerdot}(P^{\centerdot}, \sigma_{\leq n+1}X^{\centerdot})
& \to &
\sigma_{\leq n+1}X^{\centerdot} & \to &
\varDelta_{\infty}^{\centerdot}(P^{\centerdot}, \sigma_{\leq n+1}X^{\centerdot})
& \to &
\nabla_{\infty}^{\centerdot}(P^{\centerdot}, \sigma_{\leq n+1}X^{\centerdot})[1] .
\end{array}\]}\par\noindent
Then we define 
$\nabla^{\centerdot}_{\infty}(P^{\centerdot},X^{\centerdot})$ and
$\varDelta^{\centerdot}_{\infty}(P^{\centerdot},X^{\centerdot})$ by the complex of Lemma \ref{hlim1} (1)
and $\hclim \varDelta_{\infty}^{\centerdot}(P^{\centerdot}, \sigma_{\leq n}X^{\centerdot})$, 
respectively.  Moreover, we have a triangle:
\[
\nabla^{\centerdot}_{\infty}(P^{\centerdot},X^{\centerdot}) \to
X^{\centerdot} \to
\varDelta^{\centerdot}_{\infty}(P^{\centerdot},X^{\centerdot}) \to
\nabla^{\centerdot}_{\infty}(P^{\centerdot},X^{\centerdot})[1] ,
\]
because $X^{\centerdot} \cong \hclim \ \sigma_{\leq n}X^{\centerdot}$.
\end{defn}

\begin{prop} \label{unique01}
Let $A$, $B$ be projective algebras over a commutative ring $k$, 
$P^{\centerdot} \in \cat{D}(A)$ a partial tilting complex with $B \cong \End_{\cat{D}(A)}(P^{\centerdot})$,
and ${}_{B}V^{\centerdot}_{A}$ the associated bimodule complex of $P^{\centerdot}$.
For $X^{\centerdot} \in \cat{D}(A)$, we have an isomorphism of triangles in $\cat{D}(A)$:
\[\begin{CD}
j_{V !}j^{*}_{V}X^{\centerdot} @>>> X^{\centerdot} @>>> i_{V *}i^{*}_{V}X^{\centerdot}
@>>> j_{V !}j^{*}_{V}X^{\centerdot}[1] \\
@VV\wr V @| @VV\wr V @VV\wr V \\
\nabla^{\centerdot}_{\infty}(P^{\centerdot},X^{\centerdot}) @>>>
X^{\centerdot} @>>>
\varDelta^{\centerdot}_{\infty}(P^{\centerdot},X^{\centerdot}) @>>>
\nabla^{\centerdot}_{\infty}(P^{\centerdot},X^{\centerdot})[1] .
\end{CD}\]
\end{prop}

\begin{proof}
By Lemma \ref{unique00}, 
$\nabla_{\infty}^{\centerdot}(P^{\centerdot}, \sigma_{\leq n}X^{\centerdot}) \in \opn{Im}j_{V !}$
and 
$\varDelta_{\infty}^{\centerdot}(P^{\centerdot}, \sigma_{\leq n}X^{\centerdot}) \in \opn{Im}i_{V *}$.
Since $P^{\centerdot}$ is a perfect complex, $\Hom_{\cat{D}(A)}(P^{\centerdot}, -)$ commutes with
coproducts.  Then we have $\varDelta^{\centerdot}_{\infty}(P^{\centerdot},X^{\centerdot})
\in \opn{Im}i_{V *}$.
We have also $\nabla^{\centerdot}_{\infty}(P^{\centerdot},X^{\centerdot})
\in \opn{Im}j_{V !}$, because $j_{V !}$ is fully faithful and commutes with coproducts.
By Proposition \ref{t-st1}, we complete the proof.
\end{proof}

\begin{cor} \label{vardelta2}
Let $A$, $B$ be projective algebras over a commutative ring $k$, 
$P^{\centerdot} \in \cat{D}(A)$ a partial tilting complex with $B \cong \End_{\cat{D}(A)}(P^{\centerdot})$,
and ${}_{B}V^{\centerdot}_{A}$ the associated bimodule complex of $P^{\centerdot}$.
For $X^{\centerdot} \in \cat{D}(A)$, we have isomorphisms in $\cat{D}(A)$:
\[\begin{aligned}
X^{\centerdot}\lten_{A}V^{\star \centerdot}\lten_{B}V^{\centerdot}
& \cong 
\nabla^{\centerdot}_{\infty}(P^{\centerdot},X^{\centerdot}), \\
X^{\centerdot}\lten_{A}\varDelta^{\centerdot}_{A}(V^{\centerdot}) 
& \cong 
\varDelta^{\centerdot}_{\infty}(P^{\centerdot},X^{\centerdot}).
\end{aligned}\]
\end{cor}

\begin{proof}
By Theorem \ref{parttiltrecoll1} and Proposition \ref{unique01}, we complete the proof.
\end{proof}

For an idempotent $e$ of a ring $A$,
by $\Hom_{A}(eA,A) \cong Ae$, we have 
\[\begin{aligned}
j^{e}_{A !} & = - \lten_{eAe}eA : \cat{D}(eAe) \to \cat{D}(A) , \\
j_{A}^{e *} & = -\ten_{A}Ae \cong \Hom_{A}(eA,-) : \cat{D}(A) \to \cat{D}(eAe) , \\
j^{e}_{A *} & = \rhom_{eAe}(Ae,-) : \cat{D}(eAe) \to \cat{D}(A).
\end{aligned}\]
And we also get
the triangle $\xi_{e}$ in $\cat{D}(A^{\mrm{e}})$:
\[
Ae\lten_{eAe}eA \xarr{\varepsilon_{e}} A  \xarr{\eta_{e}}
\varDelta^{\centerdot}_{A}(e) \to Ae\lten_{eAe}eA[1].
\]
Throughout this paper, we identify $\cat{Mod}A/AeA$ with the full subcategory
of $\cat{Mod}A$ consisting of $A$-modules $M$ such that $\Hom_{A}(eA, M)=0$.
We denote by $\cat{D}_{A/AeA}^{*}(A)$ the full subcategory of $\cat{D}^{*}(A)$
consisting of complexes whose cohomologies are in $\cat{Mod}A/AeA$, where $*=$
nothing, $+, -, \mrm{b}$.  According to Theorem \ref{parttiltrecoll1}, we have the following.

\begin{prop} \label{idemprecoll1}
Let $A$ be a projective algebra over a commutative ring $k$, $e$ an idempotent of $A$, and let
{\small \[\begin{aligned}
j^{e}_{A !} & = - \lten_{eAe}eA : \cat{D}(eAe) \to \cat{D}(A) ,
& j_{A}^{e *} & = -\ten_{A}Ae : \cat{D}(A) \to \cat{D}(eAe) , \\
j^{e}_{A *} & = \rhom_{eAe}(Ae, -) : \cat{D}(eAe) \to \cat{D}(A) , 
& \ i^{e *}_{A} & = - \lten_{A}\varDelta^{\centerdot}_{A}(e) : \cat{D}(A) \to \cat{D}_{eAe}(A) , \\
i^{e}_{A *} & =\ \text{the embedding }: \cat{D}_{eAe}(A) \to \cat{D}(A) ,
& i^{e !}_{A} & = \rhom_{A}(\varDelta^{\centerdot}_{A}(e), - ) : \cat{D}(A) \to \cat{D}_{eAe}(A) .
\end{aligned}\]}\par\noindent
Then $\{\cat{D}_{A/AeA}(A), \cat{D}(A), \cat{D}(eAe);  
i^{e *}_{A}, i^{e}_{A *}, i^{e !}_{A}, j^{e}_{A !}, j_{A}^{e *}, j^{e}_{A *}\}$ is a recollement.
\end{prop}

\begin{rem} \label{idempbi}
According to Proposition \ref{ddotadj} and Lemma \ref{pfadj},
it is easy to see that $\{\cat{D}_{C^{\circ}\ten A/AeA}(C^{\circ}\ten A), \cat{D}(C^{\circ}\ten A), 
\cat{D}(C^{\circ}\ten eAe);  i^{e *}_{A}, i^{e}_{A *}, i^{e !}_{A}, j^{e}_{A !}, j_{A}^{e *}, j^{e}_{A *}\}$
is also a recollement for any projective $k$-algebra $C$.  
\end{rem}

\begin{cor} \label{idemprecoll2}
Let $A$ be a projective algebra over a commutative ring $k$, and $e$ an idempotent of $A$, then
the following hold.
\begin{enumerate}
\item $\varDelta^{\centerdot}_{A}(e)\lten_{A}\varDelta^{\centerdot}_{A}(e) \cong
\varDelta^{\centerdot}_{A}(e)$ in $\cat{D}(A^{\mrm{e}})$
\item $\rhom_{A}(\varDelta^{\centerdot}_{A}(e), \varDelta^{\centerdot}_{A}(e)) \cong
\varDelta^{\centerdot}_{A}(e)$ in $\cat{D}(A^{\mrm{e}})$
\item We have the following isomorphisms in $\cat{Mod}A^{\mrm{e}}$:
\[
A/AeA \cong \opn{End}_{\cat{D}(A)}(\varDelta^{\centerdot}_{A}(e))
\cong \opn{H}^{0}(\varDelta^{\centerdot}_{A}(e))  .
\]
Moreover, the first isomorphism is a ring isomorphism.
\end{enumerate}
\end{cor}

\begin{proof}
1, 2.  By Corollary \ref{vardelta1}.

3.  Applying $\Hom_{\cat{D}(A)}(-, \varDelta^{\centerdot}_{A}(e))$ to $\xi_e$,
we have an isomorphism in $\cat{Mod}A^{\mrm{e}}$:
\[
\Hom_{\cat{D}(A)}(\varDelta^{\centerdot}_{A}(e), \varDelta^{\centerdot}_{A}(e))
\cong \Hom_{\cat{D}(A)}(A, \varDelta^{\centerdot}_{A}(e)) ,
\]
because $\Hom_{\cat{D}(A)}(Ae\lten_{eAe}eA, \varDelta^{\centerdot}_{A}(e)[n]) \cong
\Hom_{\cat{D}(A)}(j^{e}_{A !}j^{e *}_{A}(A), i^{e}_{A *}i^{e !}_{A}(A)[n]) = 0$
for all $n \in \mbb{Z}$ by Proposition \ref{bilocal}, 1.  
Applying $\Hom_{\cat{D}(A)}(A, -)$ to $\xi_e$, we have an isomorphism
between exact sequences in $\cat{Mod}A^{\mrm{e}}$:
{\small\[\begin{array}{ccccc}
\Hom_{\cat{D}(A)}(A, Ae\lten_{eAe}eA) & \to & \Hom_{\cat{D}(A)}(A, A) & \to
& \Hom_{\cat{D}(A)}(A, \varDelta^{\centerdot}_{A}(e)) \to 0 \\
\downarrow\wr & & \downarrow\wr & & \downarrow\wr \\
Ae\ten_{eAe}eA & \longrightarrow & A & \longrightarrow & A/AeA \to 0 .
\end{array}\]}\par\noindent
Consider the inverse of $\Hom_{\cat{D}(A)}(\varDelta^{\centerdot}_{A}(e), \varDelta^{\centerdot}_{A}(e))
\xarr{\sim} \Hom_{\cat{D}(A)}(A, \varDelta^{\centerdot}_{A}(e))$, then
it is easy to see that $\Hom_{\cat{D}(A)}(A, A) \to \Hom_{\cat{D}(A)}(A, \varDelta^{\centerdot}_{A}(e))
\to \Hom_{\cat{D}(A)}(\varDelta^{\centerdot}_{A}(e), \varDelta^{\centerdot}_{A}(e))$ is a ring
morphism.
\end{proof}

\begin{rem} \label{idemprecoll3}
It is not hard to see that the above triangle $\xi_{e}$ also play the same role in the left
module version of Corollary \ref{idemprecoll2}.  Then we have also
\begin{enumerate}
\item $\rhom_{A^{\circ}}(\varDelta^{\centerdot}_{A}(e), \varDelta^{\centerdot}_{A}(e)) \cong
\varDelta^{\centerdot}_{A}(e)$ in $\cat{D}(A^{\mrm{e}})$
\item We have a ring isomorphism
$(A/AeA)^{\circ} \cong \opn{End}_{\cat{D}(A^{\circ})}(\varDelta^{\centerdot}_{A}(e))$.
\end{enumerate}
\end{rem}

\section{Equivalences between Recollements} \label{s3}

In this section, we study triangle equivalences between recollements induced by idempotents.

\begin{defn}
Let $\{\mcal{D}_{n}, \mcal{D}_{n}''; j_{n *}, j_{n}^{*}\}$ 
(resp., $\{\mcal{D}_{n}, \mcal{D}_{n}''; j_{n !}, j_{n}^{*}, j_{n *}\}$) be a colocalization 
(resp., a bilocalization) of $\mcal{D}_{n}$ ($n=1, 2$).
If there are triangle equivalences $F:\mcal{D}_{1} \to \mcal{D}_{2}$, $F'':\mcal{D}_{1}'' \to \mcal{D}_{2}''$
such that all squares are commutative up to ($\partial$-functorial) isomorphism in the diagram:
\[\begin{array}{cccccccc}
\ \mcal{D}_{1} & \leftrightarrows & \mcal{D}_{1}'' & & &
\ \mcal{D}_{1} & \overset{\leftarrow}{\rightleftarrows} & \mcal{D}_{1}'' \\
F \downarrow & & \ \downarrow F'' & & (\text{resp.,} & F \downarrow & & \ \downarrow F''\ ), \\
\ \mcal{D}_{2} & \leftrightarrows & \mcal{D}_{2}'' & & &
\ \mcal{D}_{2} &  \overset{\leftarrow}{\rightleftarrows} & \mcal{D}_{2}''
\end{array}\]
then we say that a colocalization $\{\mcal{D}_{1}, \mcal{D}_{1}''; j_{n *}, j_{1}^{*}\}$ 
(resp., a bilocalization $\{\mcal{D}_{1}, \mcal{D}_{1}''; \\ j_{1 !}, j_{1}^{*}, j_{1 *}\}$)
is triangle equivalent to a colocalization $\{\mcal{D}_{2}, \mcal{D}_{2}''; j_{n *}, j_{2}^{*}\}$
(resp., a bilocalization  $\{\mcal{D}_{2}, \mcal{D}_{2}''; j_{n !}, j_{2}^{*}, j_{2 *}\}$).

For recollements 
$\{\mcal{D}_{n}', \mcal{D}_{n}, \mcal{D}_{n}''; i_{n}^{*}, i_{n *},i_{n}^{!}, j_{n !}, j_{n}^{*}, j_{n *}\}$ 
($n=1,2$), if there are triangle equivalences $F':\mcal{D}_{1}' \to \mcal{D}_{2}'$, 
$F:\mcal{D}_{1} \to \mcal{D}_{2}$, $F'':\mcal{D}_{1}'' \to \mcal{D}_{2}''$
such that all squares are commutative up to ($\partial$-functorial) isomorphism in the diagram:
\[\begin{array}{ccccc}
\ \mcal{D}_{1}' & \overset{\leftarrow}{\rightleftarrows} &
\ \mcal{D}_{1} & \overset{\leftarrow}{\rightleftarrows} & \mcal{D}_{1}'' \\
F' \downarrow & & F \downarrow & & \ \downarrow F'' \\
\ \mcal{D}_{2}' & \overset{\leftarrow}{\rightleftarrows} &
\ \mcal{D}_{2} & \overset{\leftarrow}{\rightleftarrows} & \mcal{D}_{2}'' ,\\
\end{array}\]
then we say that a recollement
$\{\mcal{D}_{1}', \mcal{D}_{1}, \mcal{D}_{1}''; i_{1}^{*}, i_{1 *},i_{1}^{!}, j_{1 !}, j_{1}^{*}, j_{1 *}\}$
is triangle equivalent to a recollement
$\{\mcal{D}_{2}', \mcal{D}_{2}, \mcal{D}_{2}''; i_{2}^{*}, i_{2 *},i_{2}^{!}, j_{2 !}, j_{2}^{*}, j_{2 *}\}$.
\end{defn}

We simply write a localization $\{\mcal{D}, \mcal{D}''\}$, etc. for a localization
$\{\mcal{D}, \mcal{D}''; j^{*}, j_{*}\}$, etc. when we don't confuse them.
Parshall and Scott showed the following.
 
\begin{prop}[\cite{PS}] \label{recol}
Let $\{\mcal{D}_{n}', \mcal{D}_{n}, \mcal{D}_{n}''\}$ be recollements ($n=1,2$).  
If triangle equivalences $F:\mcal{D}_{1} \to \mcal{D}_{2}$,
$F'':\mcal{D}_{1}'' \to \mcal{D}_{2}''$ induce that a bilocalization $\{\mcal{D}_{1}, \mcal{D}_{1}''\}$
is triangle equivalent to  a bilocalization $\{\mcal{D}_{2}, \mcal{D}_{2}''\}$,
then there exists a  unique triangle equivalence $F':\mcal{D}_{1}' \to \mcal{D}_{2}'$ 
up to isomorphism such that
$F', F, F''$ induce that a recollement $\{\mcal{D}_{1}', \mcal{D}_{1}, \mcal{D}_{1}''\}$
is triangle equivalent to a  recollement $\{\mcal{D}_{2}', \mcal{D}_{2}, \mcal{D}_{2}''\}$.
\end{prop}

\begin{lem} \label{idemp1}
Let $A$ be a projective algebra over a commutative ring $k$, and $e$ an idempotent of $A$.
For $X^{\centerdot} \in \cat{D}(A)_{\mrm{perf}}$, the following are equivalent.
\begin{enumerate}
\item $X^{\centerdot} \cong P^{\centerdot}$ in $\cat{D}(A)$ 
for some $P^{\centerdot} \in \cat{K}^{\mrm{b}}(\cat{add}eA)$.
\item $j^{e}_{A !}j^{e *}_{A}(X^{\centerdot}) \cong X^{\centerdot}$ in $\cat{D}(A)$.
\item $\gamma_{X}$ is an isomorphism, where $\gamma:j^{e}_{A !}j^{e *}_{A} \to \bsym{1}_{\cat{D}(A)}$
is the adjunction arrow.
\end{enumerate}
\end{lem}

\begin{proof}
1 $\Rightarrow$ 2.
Since $j^{e}_{A !}j^{e *}_{A}(P) \cong P$ in $\cat{Mod}A$ for any $P \in \cat{add}eA$, it is trivial.

2 $\Leftrightarrow$ 3.
By Corollary \ref{t-st3}.

3 $\Rightarrow$ 1.
Let $\{Y^{\centerdot}_{i}\}_{i \in I}$ be a family  of complexes of $\cat{D}(A)$.
By Proposition \ref{pfcp}, we have isomorphisms:
\[\begin{aligned}
\coprod_{i \in I}\Hom_{\cat{D}(eAe)}(j^{e *}_{A}(X^{\centerdot}), j^{e *}_{A}(Y^{\centerdot}_{i}))
& \cong 
\coprod_{i \in I}\Hom_{\cat{D}(A)}(j^{e}_{A !}j^{e *}_{A}(X^{\centerdot}), Y^{\centerdot}_{i}) \\
& \cong
\coprod_{i \in I}\Hom_{\cat{D}(A)}(X^{\centerdot}, Y^{\centerdot}_{i}) \\
& \cong
\Hom_{\cat{D}(A)}(X^{\centerdot}, \coprod_{i \in I}Y^{\centerdot}_{i}) \\
& \cong 
\Hom_{\cat{D}(A)}(j^{e}_{A !}j^{e *}_{A}(X^{\centerdot}), \coprod_{i \in I}Y^{\centerdot}_{i}) \\
& \cong 
\Hom_{\cat{D}(eAe)}(j^{e *}_{A}(X^{\centerdot}), j^{e *}_{A}(\coprod_{i \in I}Y^{\centerdot}_{i})) \\
& \cong 
\Hom_{\cat{D}(eAe)}(j^{e *}_{A}(X^{\centerdot}), \coprod_{i \in I}j^{e *}_{A}(Y^{\centerdot}_{i})).
\end{aligned}\]
Since any complex $Z^{\centerdot}$ of $\cat{D}(eAe)$ is isomorphic to $j^{e *}_{A}(Y^{\centerdot})$ for some
$Y^{\centerdot} \in \cat{D}(A)$, by Proposition \ref{pfcp} the above isomorphisms imply
that $j^{e *}_{A}(X^{\centerdot})$ is a perfect complex of $\cat{D}(eAe)$.
Therefore, $j^{e}_{A !}j^{e *}_{A}(X^{\centerdot})$ is isomorphic to $P^{\centerdot}$
for some $P^{\centerdot} \in \cat{K}^{\mrm{b}}(\cat{add}eA)$.
\end{proof}

\begin{lem} \label{idemp2}
Let $A$, $B$ be projective algebras over a commutative ring $k$, 
and $e$, $f$ idempotents of $A$, $B$, respectively.
For $X^{\centerdot}, Y^{\centerdot} \in \cat{D}(B^{\circ}\ten A)$, 
we have an isomorphism in $\cat{D}((fBf)^{\mrm{e}})$:
\[
fB\ten _B\rhom_{A}(X^{\centerdot}, Y^{\centerdot})\ten_{B}Bf \cong
\rhom_{A}(fX^{\centerdot}, fY^{\centerdot}) .
\]
\end{lem}

\begin{proof}
First, by Proposition \ref{ddotadj}, 2, we have  isomorphisms in $\cat{D}((fBf)^{\circ}\ten B)$:
\[\begin{aligned}
fB\ten _B\rhom_{A}(X^{\centerdot}, Y^{\centerdot}) 
& \cong
\Hom_B(Bf, \rhom_{A}(X^{\centerdot}, Y^{\centerdot})) \\
& \cong
\rhom_{A}(X^{\centerdot}, \Hom_{B}(Bf, Y^{\centerdot})) \\
& \cong \rhom_{A}(X^{\centerdot}, fY^{\centerdot}).
\end{aligned}\]
Then we have  isomorphisms in $\cat{D}((fBf)^{\mrm{e}})$:
\[\begin{aligned}
fB\ten _B\rhom_{A}(X^{\centerdot}, Y^{\centerdot})\ten_{B}Bf & \cong
\rhom_{A}(X^{\centerdot}, fY^{\centerdot})\ten_{B}Bf \\
& \cong \Hom_B(fB, \rhom_{A}(X^{\centerdot}, fY^{\centerdot})) \\
& \cong \rhom_{A}(fX^{\centerdot}, fY^{\centerdot}) .
\end{aligned}\]
\end{proof}

\begin{thm} \label{recolltilt1}
Let $A, B$ be projective algebras over a commutative ring $k$, 
and $e$, $f$ idempotents of $A$, $B$, respectively.
Then the following are equivalent.
\begin{enumerate}
\item  A colocalization $\{\cat{D}(A), \cat{D}(eAe);  j^{e}_{A !}, j^{e *}_{A}\}$ is triangle equivalent to
a colocalization  $\{\cat{D}(B), \cat{D}(fBf);  j^{f}_{B !}, j^{f *}_{B}\}$.
\item There is a tilting complex $P^{\centerdot} \in \cat{K}^{\mrm{b}}(\cat{proj}A)$ such that $P^{\centerdot}=
P_{1}^{\centerdot}\oplus P_{2}^{\centerdot}$ in $\cat{K}^{\mrm{b}}(\cat{proj}A)$ satisfying \\
(a)  $B \cong \opn{End}_{\cat{D}(A)}(P^{\centerdot})$, \\
(b) under the isomorphism of (a), $f \in B$ corresponds to the canonical morphism 
$P^{\centerdot} \to P_{1}^{\centerdot} \to P^{\centerdot} \in \opn{End}_{\cat{D}(A)}(P^{\centerdot})$, \\
(c)  $P_{1}^{\centerdot} \in \cat{K}^{\mrm{b}}(\cat{add}eA)$, and 
$j^{e *}_{A}(P_{1}^{\centerdot})$ is a tilting complex for $eAe$.
\item A recollement $\{\cat{D}_{A/AeA}(A), \cat{D}(A), \cat{D}(eAe)\}$ is triangle equivalent to
a recollement $\{\cat{D}_{B/BfB}(B), \cat{D}(B), \cat{D}(fBf)\}$.
\end{enumerate}
\end{thm}

\begin{proof}
1 $\Rightarrow$ 2.
Let $G:\cat{D}(B) \to \cat{D}(A)$, $G'':\cat{D}(fBf) \to \cat{D}(eAe)$ be triangle equivalences
such that
\[\begin{array}{ccc}
\ \cat{D}(B) & \leftrightarrows & \cat{D}(fBf) \\
G \downarrow & & \ \downarrow G'' \\
\ \cat{D}(A) & \leftrightarrows & \cat{D}(eAe) \\
\end{array}\]
is commutative up to isomorphism.  Then $G(B)$ and $G''(fBf)$ are
tilting complexes for $A$ and for $eAe$ with $B \cong \opn{End}_{\cat{D}(A)}(G(B))$,
$fBf \cong \opn{End}_{\cat{D}(eAe)}(G''(B))$, respectively.  Considering $G(B)=G(fB)\oplus G((1-f)B)$,
by the above commutativity, we have isomorphisms:
\[\begin{aligned}
G(fB) & \cong Gj^{f}_{B !}(fBf) \\
& \cong j^{e}_{A !}G''(fBf) \\
& \cong j^{e}_{A !}G''j^{f *}_{B}(fB) \\
& \cong j^{e}_{A !}j^{e *}_{A}G(fB) ,
\end{aligned}\]
\[\begin{aligned}
j^{e *}_{A}G(fB) & \cong G''j^{f *}_{B}(fB) \\
& \cong G''(fBf).
\end{aligned}\]
By Proposition \ref{idemp1}, $G(fB)$ is isomorphic to a complex of $\cat{K}^{\mrm{b}}(\cat{add}eA)$,
and $j^{e *}_{A}G(fB)$ is a tilting complex for $eAe$.

2 $\Rightarrow$ 3.
Let ${}_{B}T^{\centerdot}_{A}$ be a two-sided tilting complex which is induced by $P^{\centerdot}_{A}$.
By the assumption,
$Res_{A}(fT^{\centerdot}) \cong P_{1}^{\centerdot}$ in $\cat{D}(A)$.  
By Lemma \ref{idemp1}, $\gamma_{fT}:j^{e}_{A !}j^{e *}_{A}(fT^{\centerdot}) \xarr{\sim} fT^{\centerdot}$
is an isomorphism in $\cat{D}(A)$.
By Remark \ref{idempbi}, Proposition \ref{ddotadj}, 5,
we have $fT^{\centerdot}e\lten_{eAe}eA \cong fT^{\centerdot}$ in $\cat{D}((fBf)^{\circ}\ten A)$.
By Proposition \ref{tilt2}, Lemma \ref{idemp2}, we have isomorphisms in $\cat{D}((fBf)^{\mrm{e}})$:
\[\begin{aligned}
fBf & \cong \rhom_{A}(fT^{\centerdot}, fT^{\centerdot}) \\
& \cong \rhom_{A}(fT^{\centerdot}e\lten_{eAe}eA, fT^{\centerdot}e\lten_{eAe}eA) \\
& \cong \rhom_{A}(fT^{\centerdot}e, fT^{\centerdot}e\lten_{eAe}eAe) \\
& \cong \rhom_{eAe}(fT^{\centerdot}e, fT^{\centerdot}e).
\end{aligned}\]
By taking cohomology, we have 
\[\begin{aligned}
fBf & \cong \Hom_{\cat{D}(eAe)}(fT^{\centerdot}e, fT^{\centerdot}e).
\end{aligned}\]
By the assumption, $fT^{\centerdot}e \cong j^{e *}_{A}(fT^{\centerdot}) 
\cong j^{e *}_{A}(P_{1}^{\centerdot})$ is a tilting complex
for $eAe$.  Since it is easy to see the above isomorphism is induced by the left multiplication,
by \cite{Rd2} Lemma 3.2, \cite{Ke} Theorem, $fT^{\centerdot}e$ 
is a two-sided tilting complex in $\cat{D}((fBf)^{\circ}\ten eAe)$.
Let
{\small \[\begin{aligned}
F & = \rhom_{A}(T^{\centerdot}, - ) : \cat{D}(B^{\circ}\ten A) \to \cat{D}(B^{\circ}\ten B) , \\
F'' & = \rhom_{eAe}(fT^{\centerdot}e, - ) : \cat{D}(B^{\circ}\ten eAe) \to \cat{D}(B^{\circ}\ten fBf) , \\
G  & = - \lten_{B}T^{\centerdot} : \cat{D}(B^{\circ}\ten B) \to \cat{D}(B^{\circ}\ten A) , \\
G'' & = - \lten_{fBf}fT^{\centerdot}e : \cat{D}(B^{\circ}\ten eAe) \to \cat{D}(B^{\circ}\ten fBf) .
\end{aligned}\]}\par\noindent
Using the same symbols, consider a triangle equivalence between colocalizations 
$\{\cat{D}(B^{\circ}\ten A), \cat{D}(B^{\circ}\ten eAe); j^{e}_{A !}, j_{A}^{e *} \}$
and $\{\cat{D}(B^{\circ}\ten B), \cat{D}(B^{\circ}\ten fBf); j^{f}_{B !}, j_{B}^{f *} \}$.
And we use the same symbols
{\small \[\begin{aligned}
F & = \rhom_{A}(T^{\centerdot}, - ) : \cat{D}(A) \to \cat{D}(B) ,
& \ F'' & = \rhom_{eAe}(fT^{\centerdot}e, - ) : \cat{D}(eAe) \to \cat{D}(fBf) , \\
G  & = - \lten_{B}T^{\centerdot} : \cat{D}(B) \to \cat{D}(A) ,
& G'' & = - \lten_{fBf}fT^{\centerdot}e : \cat{D}(eAe) \to \cat{D}(fBf) .
\end{aligned}\]}\par\noindent
For any $X^{\centerdot} \in \cat{D}(B^{\circ}\ten A)$ (resp., $X^{\centerdot} \in \cat{D}(A)$),
by Proposition \ref{ddotadj}, 3,
we have isomorphisms in $\cat{D}(B^{\circ}\ten fBf)$ (resp., $\cat{D}(fBf)$):
\[\begin{aligned}
j^{f *}_{B}F(X^{\centerdot}) & \cong \rhom_{B}(fB, \rhom_{A}(T^{\centerdot}, X^{\centerdot})) \\
& \cong \rhom_{A}(fT^{\centerdot}, X^{\centerdot}) \\
& \cong \rhom_{A}(j^{e}_{A !}j^{e *}_{A}(fT^{\centerdot}), X^{\centerdot}) \\
& \cong \rhom_{eAe}(j^{e *}_{A}(fT^{\centerdot}), j^{e *}_{A}(X^{\centerdot})) \\
& \cong F''j^{e *}_{A}(X^{\centerdot}) .
\end{aligned}\]
Since $G$, $G''$ are quasi-inverses of $F$, $F''$, respectively, for $B \in \cat{D}(B^{\circ}\ten B)$
we have isomorphisms in $\cat{D}(B^{\circ}\ten eAe)$:
\[\begin{aligned}
T^{\centerdot}e & \cong j_{A}^{e *}G(B) \\
& \cong G''j_{B}^{f *}(B) \\
& \cong Bf\lten_{fBf}fT^{\centerdot}e
\end{aligned}\]
Therefore, for any $Y^{\centerdot} \in \cat{D}(eAe)$, we have isomorphisms in $\cat{D}(B)$:
\[\begin{aligned}
j^{f }_{B *}F''(Y^{\centerdot}) & \cong \rhom_{fBf}(Bf , \rhom_{eAe}(fT^{\centerdot}e, Y^{\centerdot})) \\
& \cong \rhom_{B}(Bf\lten_{fBf}fT^{\centerdot}e, Y^{\centerdot}) \\
& \cong \rhom_{B}(T^{\centerdot}e, Y^{\centerdot}) \\
& \cong \rhom_{B}(j^{e *}_{A}(T^{\centerdot}), Y^{\centerdot}) \\
& \cong \rhom_{B}(T^{\centerdot}, j^{e}_{A *}(Y^{\centerdot})) \\
& \cong Fj^{e}_{A *}(Y^{\centerdot}).
\end{aligned}\]
For any $Z^{\centerdot} \in \cat{D}(fBf)$, we have isomorphisms in $\cat{D}(A)$:
\[\begin{aligned}
j^{e}_{A !}G''(Z^{\centerdot}) & = Z^{\centerdot}\lten_{fBf}fT^{\centerdot}e\lten_{eAe}eA \\
& \cong Z^{\centerdot}\lten_{fBf}fT^{\centerdot} \\
& \cong Z^{\centerdot}\lten_{fBf}fB\ten_{B}T^{\centerdot} \\
& \cong G''j^{f}_{B !}(Z^{\centerdot}).
\end{aligned}\]
Since $F$, $F''$ are quasi-inverses of $G$, $G''$, respectively, we have
$j^{f}_{B !}F'' \cong Fj^{e}_{A !}$.
By Proposition \ref{recol}, we have the statement.

3 $\Rightarrow$ 1.
It is trivial.
\end{proof}

\begin{defn}
Let $A$ be a projective algebra over a commutative ring $k$, and $e$ an idempotent of $A$.
We call a tilting complex $P^{\centerdot} \in \cat{K}^{\mrm{b}}(\cat{projA})$ a recollement tilting complex 
related to an idempotent $e$ of $A$ if $P^{\centerdot}$ satisfies 
the condition of Theorem \ref{recolltilt1}, 2.  In this case, we call an idempotent $f \in B$ 
an idempotent corresponding to $e$.
\end{defn}

We see the following symmetric properties of a two-sided tilting complex which is induced by 
a recollement tilting complex.  We will call the following two-sided tilting complex 
a {\it two-sided recollement tilting complex ${}_{B}T^{\centerdot}_{A}$ related to idempotents $e \in A$, $f \in B$}.

\begin{cor} \label{recolltilt2}
Let $A, B$ be projective algebras over a commutative ring $k$, and $e$, $f$ idempotents of $A$, $B$, respectively.
Let ${}_{B}T^{\centerdot}_{A}$ be a two-sided tilting complex
such that
\begin{enumerate}
\rmitem{a} $fT^{\centerdot}e \in \cat{D}((fBf)^{\circ}\ten eAe)$ is a two-sided tilting complex,
\rmitem{b} $fT^{\centerdot}e\lten_{eAe}eA \cong fT^{\centerdot}$ in $\cat{D}((fBf)^{\circ}\ten A)$.
\end{enumerate}
Then the following hold.
\begin{enumerate}
\item  $Bf\lten_{fBf}fT^{\centerdot}e \cong T^{\centerdot}e$ in $\cat{D}(B^{\circ}\ten eAe)$.
\item  $eT^{\vee \centerdot}f$ is the inverse of $fT^{\centerdot}e$, where $T^{\vee \centerdot}$ is 
the inverse of $T^{\centerdot}$.
\item  $Ae\lten_{eAe}eT^{\vee \centerdot}f \cong T^{\vee \centerdot}f$ in $\cat{D}(A^{\circ}\ten fBf)$.
\item  $eT^{\vee \centerdot}f\lten_{fBf}fB \cong eT^{\vee \centerdot}$ in $\cat{D}((eAe)^{\circ}\ten B)$.
\end{enumerate}
\end{cor}

\begin{proof}
Here we use the same symbols in the proof 2 $\Rightarrow$ 3 of Theorem \ref{recolltilt1}.
It is easy to see that $F$ and $F''$ induce a triangle equivalence between bilocalizations 
$\{\cat{D}(B^{\circ}\ten A), \cat{D}(B^{\circ}\ten eAe); j^{e}_{A !}, j_{A}^{e *}, j^{e}_{A *} \}$
and $\{\cat{D}(B^{\circ}\ten B), \cat{D}(B^{\circ}\ten fBf); j^{f}_{B !}, j_{B}^{f *}, j^{f}_{B *} \}$.
By the proof of Theorem \ref{recolltilt1}, we get the statement 1, and 
$j^{f *}_{B}F \cong F''j^{e *}_{A}$, $j^{f}_{B !}F'' \cong Fj^{e}_{A !}$ and 
$j^{f}_{B *}F'' \cong Fj^{e}_{A *}$.
Then we have  isomorphisms $j^{f *}_{B}Fj^{e}_{A !} \cong F''j^{e *}_{A}j^{e}_{A !} \cong F''$.
Since $-\lten_{A}T^{\vee \centerdot}_{B} \cong F$, we have isomorphisms
$eT^{\vee \centerdot}f \cong \rhom_{eAe}(fT^{\centerdot}e, eAe)$ in $\cat{D}((eAe)^{\circ}\ten fBf)$,
and $- \lten_{eAe}eT^{\vee \centerdot}f \cong F''$.
This means that $eT^{\vee \centerdot}f$ is the inverse of a two-sided tilting complex $fT^{\centerdot}e$.
Similarly, $j^{f *}_{B}F \cong F''j^{e *}_{A}$ and $j^{f}_{B !}F'' \cong Fj^{e}_{A !}$ imply
the statements 3 and 4, respectively.
\end{proof}

\begin{cor} \label{recolltilt3}
Let $A, B$ be projective algebras over a commutative ring $k$, and $e$, $f$ idempotents of $A$, $B$, respectively.
For a two-sided recollement tilting complex ${}_{B}T^{\centerdot}_{A}$ related to 
idempotents $e$, $f$, we have an isomorphism between triangles
$T^{\centerdot}\lten_{A}\xi_{e}$ and $\xi_{f}\lten_{B}T^{\centerdot}$ in 
$\cat{D}(B^{\circ}\ten A)$:
\[\begin{CD}
T^{\centerdot}e\lten_{eAe}eA @>>> T^{\centerdot} @>>> 
T^{\centerdot}\lten_{A}\varDelta_{A}^{\centerdot}(e) @>>> 
T^{\centerdot}e\lten_{eAe}eA[1] \\
@VV\wr V @| @VV\wr V @VV\wr V \\
Bf\lten_{fBf}fT^{\centerdot} @>>> T^{\centerdot} @>>> 
\varDelta_{B}^{\centerdot}(f)\lten_{B}T^{\centerdot} @>>> 
Bf\lten_{fBf}fT^{\centerdot}[1] .
\end{CD}\]
\end{cor}

\begin{proof}
According to Proposition \ref{recol}, for the triangle equivalence between colocalizations 
in the proof of  Corollary \ref{recolltilt2} there exists 
$F':\cat{D}_{B^{\circ}\ten B/BfB}(B^{\circ}\ten B) \to 
\cat{D}_{B^{\circ}\ten A/AeA}(B^{\circ}\ten A)$ such that
a recollement
\[
\{\cat{D}_{B^{\circ}\ten B/BfB}(B^{\circ}\ten B), \cat{D}(B^{\circ}\ten B), 
\cat{D}(B^{\circ}\ten fBf); i^{f *}_{B}, i^{f}_{B *}, i^{f !}_{B}, 
j^{f}_{B !}, j^{f *}_{B}, j^{f}_{B *}\}
\]
is triangle equivalent to a recollement
\[
\{\cat{D}_{B^{\circ}\ten A/AeA}(B^{\circ}\ten A), \cat{D}(B^{\circ}\ten A),
\cat{D}(B^{\circ}\ten eAe); i^{e *}_{A}, i^{e}_{A *}, i^{e !}_{A}, 
j^{e}_{A !}, j^{e *}_{A}, j^{e}_{A *}\}.
\]
By Proposition \ref{ddotadj}, Lemma \ref{pfadj}, 
the triangle $T^{\centerdot}\lten_{A}\xi_{e}$ is isomorphic to the following triangle 
in $\cat{D}(B^{\circ}\ten A)$:
\[
j^{e}_{A !}j^{e *}_{A}(T^{\centerdot}) \to T^{\centerdot} \to
i^{e}_{A *}i^{e *}_{A}(T^{\centerdot}) \to j^{e}_{A !}j^{e *}_{A}(T^{\centerdot})[1] .
\]
On the other hand, the triangle $\xi_{f}\lten_{B}T^{\centerdot}$ is isomorphic to
the following triangle in $\cat{D}(B^{\circ}\ten A)$:
\[
Fj^{f}_{B !}j^{f *}_{B}(B) \to F(B) \to
Fi^{f}_{B *}i^{f *}_{B}(B) \to Fj^{f}_{B !}j^{f *}_{B}(B)[1] .
\]
Since
$F(B) \cong T^{\centerdot}$,
$Fj^{f}_{B !}j^{f *}_{B}(B) \cong j^{e}_{A !}F''j^{f *}_{B}(B) \cong j^{e}_{A !}j^{e *}_{A}F(B)$,
$Fi^{f}_{B *}i^{f *}_{B}(B) \cong i^{e}_{A *}F'i^{f *}_{B}(B) \cong i^{e}_{A *}i^{e *}_{A}F(B)$,
by Proposition \ref{t-st1}, we complete the proof.
\end{proof}

\begin{cor} \label{recolltilt4}
Let $A, B$ be projective algebras over a commutative ring $k$, and $e$, $f$ idempotents of $A$, $B$, respectively.
For a two-sided recollement tilting complex ${}_{B}T^{\centerdot}_{A}$ related to 
idempotents $e$, $f$, the following hold.
\begin{enumerate}
\item $T^{\centerdot}\lten_{A}\varDelta_{A}^{\centerdot}(e) \cong
\varDelta_{B}^{\centerdot}(f)\lten_{B}T^{\centerdot}$ in $\cat{D}(B^{\circ}\ten A)$.
\item $\varDelta_{A}^{\centerdot}(e)\lten_{A}T^{\vee \centerdot} \cong
T^{\vee \centerdot}\lten_{B}\varDelta_{B}^{\centerdot}(f)$ in $\cat{D}(A^{\circ}\ten B)$.
\end{enumerate}
\end{cor}

\begin{proof}
1.  By Corollary \ref{recolltilt3}.

2.  We have isomorphisms in $\cat{D}(A^{\circ}\ten B)$:
\[\begin{aligned}
\varDelta_{A}^{\centerdot}(e)\lten_{A}T^{\vee \centerdot}
& \cong
T^{\vee \centerdot}\lten_{B}T^{\centerdot}\lten_{A}\varDelta_{A}^{\centerdot}(e)\lten_{A}T^{\vee \centerdot} \\
& \cong
T^{\vee \centerdot}\lten_{B}\varDelta_{B}^{\centerdot}(f)\lten_{B}T^{\centerdot}\lten_{A}T^{\vee \centerdot} \\
& \cong
T^{\vee \centerdot}\lten_{B}\varDelta_{B}^{\centerdot}(f) .
\end{aligned}\]
\end{proof}

\begin{defn}
Let $A, B$ be projective algebras over a commutative ring $k$, and $e$, $f$ idempotents of $A$, $B$, respectively.
For a two-sided recollement tilting complex ${}_{B}T^{\centerdot}_{A}$ related to 
idempotents $e$, $f$, we define
\[\begin{array}{cc}
\varDelta^{\centerdot}_{T} = T^{\centerdot}\lten_{A}\varDelta_{A}^{\centerdot}(e) 
\in \cat{D}(B^{\circ}\ten A), &
\varDelta^{\vee \centerdot}_{T} = \varDelta_{A}^{\centerdot}(e)\lten_{A}T^{\vee \centerdot} 
\in \cat{D}(A^{\circ}\ten B).
\end{array}\]
\end{defn}

\begin{prop} \label{recolltilt5}
Let $A, B$ be projective algebras over a commutative ring $k$, and $e$, $f$ idempotents of $A$, $B$, respectively.
For a two-sided recollement tilting complex ${}_{B}T^{\centerdot}_{A}$ related to 
idempotents $e$, $f$, let
\[\begin{aligned}
F' & =\rhom_{A}(\varDelta^{\centerdot}_{T}, -):\cat{D}_{A/AeA}(A) 
\to \cat{D}_{B/BfB}(B) , \\
F & =\rhom_{A}(T^{\centerdot}, -): \cat{D}(A) \to \cat{D}(B) , \\
F'' & =\rhom_{eAe}(fT^{\centerdot}e, -) : \cat{D}(eAe) \to \cat{D}(fBf) . 
\end{aligned}\]
Then the following hold.
\begin{enumerate}
\item We have an isomorphism $F' \cong -\lten_{A}\varDelta^{\vee \centerdot}_{T}$.
\item A quasi-inverse $G'$ of $F'$ is isomorphic to $\rhom_{B}(\varDelta^{\vee \centerdot}_{T}, -)
\cong -\lten_{B}\varDelta^{\centerdot}_{T}$.
\item $F'$, $F$, $F''$ induce that 
a recollement $\{\cat{D}_{A/AeA}(A), \cat{D}(A), \cat{D}(eAe)\}$ is triangle equivalent to
a recollement $\{\cat{D}_{B/BfB}(B), \cat{D}(B), \cat{D}(fBf)\}$.
\end{enumerate}
\end{prop}

\begin{proof}
According to Proposition \ref{recol}, $F'$ exists and satisfies
$F' \cong i^{f *}_{B}Fi^{e}_{A *} \cong i^{f !}_{B}Fi^{e}_{A *}$.
By Proposition \ref{idemprecoll1}, we have isomorphisms
\[\begin{aligned}
i^{f *}_{B}Fi^{e}_{A *} & \cong 
\rhom_{A}(T^{\centerdot}, -)\lten_{B}\varDelta^{\centerdot}_{B}(f) \\
& \cong
-\lten_{A}T^{\vee \centerdot}\lten_{B}\varDelta^{\centerdot}_{B}(f), \\
i^{f !}_{B}Fi^{e}_{A *} & \cong 
\rhom_{B}(\varDelta^{\centerdot}_{B}(f), \rhom_{A}(T^{\centerdot}, -)) \\
& \cong 
\rhom_{A}(\varDelta^{\centerdot}_{B}(f)\lten_{A}T^{\centerdot}, -) .
\end{aligned}\]
Let $G=\rhom_{B}(T^{\vee \centerdot}, -)$.
Since $G' \cong i^{e *}_{A}Gi^{f}_{B *} \cong i^{e !}_{B}Gi^{f}_{B *}$,
we have isomorphisms
\[\begin{aligned}
i^{e *}_{A}Gi^{f}_{B *} & \cong 
\rhom_{B}(T^{\vee \centerdot}, -)\lten_{A}\varDelta^{\centerdot}_{A}(e) \\
& \cong
-\lten_{B}T^{\centerdot}\lten_{A}\varDelta^{\centerdot}_{A}(e), \\
i^{e !}_{A}Gi^{f}_{B *} & \cong 
\rhom_{A}(\varDelta^{\centerdot}_{A}(e), \rhom_{B}(T^{\vee \centerdot}, -)) \\
& \cong 
\rhom_{B}(\varDelta^{\centerdot}_{A}(e)\lten_{A}T^{\vee \centerdot}, -) .
\end{aligned}\]
By Corollary \ref{recolltilt4}, we complete the proof.
\end{proof}

\begin{cor} \label{bounded}
Under the condition of Proposition \ref{recolltilt5}, the following hold.
\begin{enumerate}
\item $Res_{A}\varDelta^{\centerdot}_{T}$ is a compact object in $\cat{D}_{A/AeA}(A)$.
\item $Res_{B^{\circ}}\varDelta^{\centerdot}_{T}$ is a compact object in 
$\cat{D}_{(B/BfB)^{\circ}}(B^{\circ})$.
\item $\rhom_{A}(\varDelta^{\centerdot}_{T}, -):\cat{D}_{A/AeA}^{*}(A) \xarr{\sim} \cat{D}_{B/BfB}^{*}(B)$
is a triangle equivalence, where $*=$ nothing, $+, -, \mrm{b}$.
\end{enumerate}
\end{cor}

\begin{proof}
1, 2.
By Corollary \ref{cpt1}, it is trivial.

3.
Since for any $X^{\centerdot} \in \cat{D}_{A/AeA}(A)$ we have isomorphisms in $\cat{D}_{B/BfB}(B)$:
\[\begin{aligned}
F'(X^{\centerdot}) & =
\rhom_{A}(\varDelta^{\centerdot}_{T}, X^{\centerdot}) \\
& =
\rhom_{A}(T^{\centerdot}\lten_{A}\varDelta_{A}^{\centerdot}(e), X^{\centerdot}) \\
& \cong
\rhom_{A}(T^{\centerdot}, \rhom_{A}(\varDelta_{A}^{\centerdot}(e), X^{\centerdot})) \\
& \cong
\rhom_{A}(T^{\centerdot}, X^{\centerdot}) , 
\end{aligned}\]
we have $\opn{Im}F'|_{\cat{D}_{A/AeA}^{*}(A)} \subset \cat{D}_{B/BfB}^{*}(B)$,
where $*=$ nothing, $+, -, \mrm{b}$.
Let $G'=\rhom_{B}(\varDelta^{\vee \centerdot}_{T}, -)$, then
we have also $\opn{Im}G'|_{\cat{D}_{B/BfB}^{*}(B)} \subset \cat{D}_{A/AeA}^{*}(A)$,
where $*=$ nothing, $+, -, \mrm{b}$.
Since $G'$ is a quasi-inverse of $F'$, we complete the proof.
\end{proof}

\begin{prop} \label{recolltilt6}
Let $A, B$ be projective algebras over a commutative ring $k$, and $e$, $f$ idempotents of $A$, $B$, respectively.
For a two-sided recollement tilting complex ${}_{B}T^{\centerdot}_{A}$ related to 
idempotents $e$, $f$, the following hold.
\begin{enumerate}
\item $\rhom_{A}(\varDelta^{\centerdot}_{T}, \varDelta^{\centerdot}_{T}) \cong
\varDelta^{\centerdot}_{T}\lten_{A}\varDelta^{\vee \centerdot}_{T} 
\cong \varDelta^{\centerdot}_{B}(f)$ in $\cat{D}(B^{\mrm{e}})$.
\item $\rhom_{B^{\circ}}(\varDelta^{\centerdot}_{T}, \varDelta^{\centerdot}_{T}) \cong
\varDelta^{\vee \centerdot}_{T}\lten_{B}\varDelta^{\centerdot}_{T} 
\cong \varDelta^{\centerdot}_{A}(e)$ in $\cat{D}(A^{\mrm{e}})$.
\item We have a ring isomorphism $\End_{\cat{D}(A)}(\varDelta^{\centerdot}_{T}) \cong B/BfB$.
\item We have a ring isomorphism $\End_{\cat{D}(B^{\circ})}(\varDelta^{\centerdot}_{T}) \cong (A/AeA)^{\circ}$.
\end{enumerate}
\end{prop}

\begin{proof}
1.  By Corollaries \ref{idemprecoll2}, \ref{recolltilt4}, Proposition \ref{recolltilt5}, 
we have isomorphisms in $\cat{D}(B^{\mrm{e}})$:
\[\begin{aligned}
\rhom_{A}(\varDelta^{\centerdot}_{T}, \varDelta^{\centerdot}_{T})
& \cong
\varDelta^{\centerdot}_{T}\lten_{A}\varDelta^{\vee \centerdot}_{T} \\
& \cong
\varDelta_{B}^{\centerdot}(f)\lten_{B}T^{\centerdot}\lten_{A}
T^{\vee \centerdot}\lten_{B}\varDelta_{B}^{\centerdot}(f) \\
& \cong
\varDelta_{B}^{\centerdot}(f)\lten_{B}\varDelta_{B}^{\centerdot}(f) \\
& \cong
\varDelta_{B}^{\centerdot}(f) .
\end{aligned}\]

2.  By Remark \ref{idemprecoll3}, Corollary \ref{idemprecoll2}, we have isomorphisms in $\cat{D}(A^{\mrm{e}})$:
\[\begin{aligned}
\rhom_{B^{\circ}}(\varDelta^{\centerdot}_{T}, \varDelta^{\centerdot}_{T})  & =
\rhom_{B^{\circ}}(T^{\centerdot}\lten_{A}\varDelta_{A}^{\centerdot}(e), 
T^{\centerdot}\lten_{A}\varDelta_{A}^{\centerdot}(e)) \\
& \cong
\rhom_{A^{\circ}}(\varDelta_{A}^{\centerdot}(e), 
\rhom_{B^{\circ}}(T^{\centerdot}, T^{\centerdot}\lten_{A}\varDelta_{A}^{\centerdot}(e))) \\
& \cong
\rhom_{A^{\circ}}(\varDelta_{A}^{\centerdot}(e), \varDelta_{A}^{\centerdot}(e)) \\
& \cong
\varDelta_{A}^{\centerdot}(e) ,
\end{aligned}\]
and have isomorphisms in $\cat{D}(A^{\mrm{e}})$:
\[\begin{aligned}
\varDelta^{\vee \centerdot}_{T}\lten_{B}\varDelta^{\centerdot}_{T} 
& \cong
\varDelta_{A}^{\centerdot}(e)\lten_{A}T^{\vee \centerdot}\lten_{B}
T^{\centerdot}\lten_{A}\varDelta_{A}^{\centerdot}(e) \\
& \cong
\varDelta_{A}^{\centerdot}(e)\lten_{A}\varDelta_{A}^{\centerdot}(e) \\
& \cong
\varDelta_{A}^{\centerdot}(e) .
\end{aligned}\]

3.  By Corollaries \ref{idemprecoll2}, \ref{recolltilt4}, we have ring isomorphisms:
\[\begin{aligned}
\End_{\cat{D}(A)}(\varDelta^{\centerdot}_{T}) & \cong
\End_{\cat{D}(B)}(\varDelta^{\centerdot}_{T}\lten_{A}T^{\vee \centerdot}) \\
& \cong
\End_{\cat{D}(B)}(\varDelta_{B}^{\centerdot}(f)\lten_{B}T^{\centerdot}\lten_{A}T^{\vee \centerdot}) \\
& \cong
\End_{\cat{D}(B)}(\varDelta_{B}^{\centerdot}(f)) \\
& \cong B/BfB .
\end{aligned}\]

4. By taking cohomology of the isomorphism of 2, we have the statement by Remark \ref{idemprecoll3}.
\end{proof}

We give some tilting complexes satisfying the following proposition in Section \ref{s4}.

\begin{prop} \label{Moritaeqv1}
Let $A$,  $B$ be projective algebras over a commutative ring $k$, $e$ an idempotent of $A$,
$P^{\centerdot}$ a recollement tilting complex related to $e$, and $B\cong\End_{\cat{D}(A)}(P^{\centerdot})$.
If $P^{\centerdot}\lten_{A}\varDelta_{A}^{\centerdot}(e) \cong \varDelta_{A}^{\centerdot}(e)$
in $\cat{D}(A)$, then the following hold.
\begin{enumerate}
\item $A/AeA \cong B/BfB$ as a ring, where $f$ is an idempotent of $B$ corresponding to $e$.
\item The standard equivalence $\rhom_{A}(T^{\centerdot}, -):\cat{D}(A) \to \cat{D}(B)$ induces
an equivalence $R^{0}\dhom_{A}(T^{\centerdot}, -)|_{\cat{Mod}A/AeA}:
\cat{Mod}A/AeA \to \cat{Mod}B/BfB$, where
${}_{B}T^{\centerdot}_{A}$ is the associated two-sided tilting complex of $P^{\centerdot}$.
\end{enumerate}
\end{prop}

\begin{proof}
1.  By the assumption, we have an isomorphism
$Res_{A}\varDelta^{\centerdot}_{T} \cong Res_{A}\varDelta_{A}^{\centerdot}(e)$
in $\cat{D}(A)$.  By Corollary \ref{idemprecoll2}, Proposition \ref{recolltilt6}, we have the statement.

2.  Let $\cat{D}^{0}_{A/AeA}(A)$ (resp., $\cat{D}^{0}_{B/BfB}(B)$) be the full subcategory
of $\cat{D}_{A/AeA}(A)$ (resp., $\cat{D}_{B/BfB}(B)$) consisting of complexes
$X^{\centerdot}$ with $\opn{H}^{i}(X^{\centerdot})=O$ for $i \not=0$.
This category is equivalent to $\cat{Mod}A/AeA$ (res., $\cat{Mod}B/BfB$).
By Corollary \ref{recolltilt4}, we have isomorphisms in $\cat{D}(B)$:
\[\begin{aligned}
\varDelta^{\vee \centerdot}_{T} & \cong 
\varDelta_{A}^{\centerdot}(e)\lten_{A}T^{\vee \centerdot} \\
& \cong 
T^{\centerdot}\lten_{A}\varDelta_{A}^{\centerdot}(e)\lten_{A}T^{\vee \centerdot} \\
& \cong 
\varDelta_{B}^{\centerdot}(f)\lten_{B}T^{\centerdot}\lten_{A}T^{\vee \centerdot} \\
& \cong \varDelta_{B}^{\centerdot}(f) .
\end{aligned}\]
Define 
\[\begin{aligned}
F' &=\rhom_{A}(\varDelta^{\centerdot}_{T}, -):
\cat{D}_{A/AeA}(A) \to \cat{D}_{B/BfB}(B) , \\
G' &=\rhom_{A}(\varDelta^{\vee \centerdot}_{T}, -):
\cat{D}_{B/BfB}(B) \to \cat{D}_{A/AeA}(A) ,
\end{aligned}\]
then they induce an equivalence between $\cat{D}_{A/AeA}(A)$ and $\cat{D}_{B/BfB}(B)$,
by Proposition \ref{recolltilt5}.
For any $X \in \cat{Mod}A/AeA$, we have isomorphisms in $\cat{D}(k)$:
\[\begin{aligned}
Res_{k}\rhom_{A}(\varDelta^{\centerdot}_{T}, X)
& \cong 
Res_{k}\rhom_{A}(\varDelta^{\centerdot}_{A}(e), X) \\
& \cong X .
\end{aligned}\]
This means that $\opn{Im}F'|_{\cat{Mod}A/AeA}$ is contained in $\cat{D}^{0}_{B/BfB}(B)$.
Similarly since we have isomorphisms in $\cat{D}(k)$:
\[\begin{aligned}
Res_{k}\rhom_{B}(\varDelta^{\vee \centerdot}_{T}, Y)
& \cong 
Res_{k}\rhom_{B}(\varDelta^{\centerdot}_{B}(f), Y) \\
& \cong Y ,
\end{aligned}\]
for any $Y \in \cat{Mod}B/BfB$, $\opn{Im}G'|_{\cat{Mod}B/BfB}$ is contained in $\cat{D}^{0}_{A/AeA}(A)$.
Therefore $F'$ and $G'$ induce an equivalence between $\cat{D}^{0}_{A/AeA}(A)$
and $\cat{D}^{0}_{B/BfB}(B)$.
Since we have isomorphisms in $\cat{D}(B)$:
\[\begin{aligned}
\rhom_{A}(T^{\centerdot}, X) & \cong
\rhom_{A}(T^{\centerdot}, i^{e}_{A *}(X)) \\
& \cong
i^{f}_{B *}\rhom_{A}(\varDelta^{\centerdot}_{T}, X)
\end{aligned}\]
for any $X \in \cat{Mod}A/AeA$, we complete the proof.
\end{proof}

\section{Tilting Complexes over symmetric Algebras} \label{s4}

Throughout this section, $A$ is a finite dimensional algebra over a field $k$,
and $D=\Hom_{k}(-, k)$.
$A$ is called a symmetric $k$-algebra if $A \cong DA$ as $A$-bimodules.
In the case of symmetric algebras, the following basic property has been seen in \cite{Rd3}.

\begin{lem} \label{symm}
Let $A$ be a symmetric algebra over a field $k$, and $P^{\centerdot} \in \cat{K}^{\mrm{b}}(\cat{proj}A)$.
For a bounded complex $X^{\centerdot}$ of finitely generated right $A$-modules, we have
an isomorphism:
\[
\dhom_{A}(P^{\centerdot}, X^{\centerdot}) \cong
D\dhom_{A}(X^{\centerdot}, P^{\centerdot}) .
\]
In particular we have an isomorphism:
\[
\Hom_{\cat{K}(A)}(P^{\centerdot}, X^{\centerdot}[n]) \cong
D\Hom_{\cat{K}(A)}(X^{\centerdot}, P^{\centerdot}[-n])
\]
for any $n \in \mbb{Z}$.
\end{lem}

\begin{defn} \label{cpxlg}
For a complex $X^{\centerdot}$, we denote
$l(X^{\centerdot})=\max \{n \ | \ \opn{H}^{n}(X^{\centerdot}) \not= 0\}
-\min \{n \ | \ \opn{H}^{n}(X^{\centerdot}) \not= 0\}+1$.
We call $l(X^{\centerdot})$ the length of a complex $X^{\centerdot}$.
\end{defn}

We redefine precisely Definition \ref{constadj} for constructing tilting complexes. 

\begin{defn} \label{consttilt}
Let $A$ be a finite dimensional algebra over a field $k$, $M$ a finitely generated $A$-module, 
and $P^{\centerdot}: P^{s-r} \to \ldots \to P^{s-1} \to P^{s} \in \cat{K}^{\mrm{b}}(\cat{proj}A)$
a partial tilting complex of length $r+1$.
For an integer $n \geq 0$, by induction, we construct a family  
$\{\varDelta^{\centerdot}_{n}(P^{\centerdot},M)\}_{n \geq 0}$ of 
complexes as follows.

Let $\varDelta^{\centerdot}_{0}(P^{\centerdot}, M) = M$.  For $n \geq 1$, by induction we construct 
a triangle $\zeta_{n}(P^{\centerdot}, M)$:
\[
P^{\centerdot}_{n}[n+s-r-1] \xarr{g_n} \varDelta^{\centerdot}_{n-1}(P^{\centerdot},M)
\xarr{h_n} \varDelta^{\centerdot}_{n}(P^{\centerdot},M) \arr P^{\centerdot}_{n}[n+s-r]
\]
as follows.  If $\opn{Hom}_{\cat{K}(A)}(P^{\centerdot}, 
\varDelta^{\centerdot}_{n-1}(P^{\centerdot},M)[r-s-n+1]) = 0$,  then we set $P^{\centerdot}_{n} = O$.  
Otherwise, we take $P^{\centerdot}_{n} \in \cat{add}P^{\centerdot}$ and a morphism 
$g'_n : P^{\centerdot}_{n}  \arr \varDelta^{\centerdot}_{n-1}(P^{\centerdot},M)[r-s-n+1]$ such 
that $\opn{Hom}_{\cat{K}(A)}(P^{\centerdot},  g'_n)$ is a projective cover as 
$\End_{\cat{D}(A)}(P^{\centerdot})$-modules, and $g_n = g'_n[n+s-r-1]$.
Moreover, $\Theta_{n}^{\centerdot}(P^{\centerdot}, M) = 
\varDelta^{\centerdot}_{n}(P^{\centerdot},M)\oplus P^{\centerdot}[n+s-r]$.
\end{defn}

By the construction, we have the following properties.

\begin{lem} \label{exttilt0}
For $\{\varDelta^{\centerdot}_{n}(P^{\centerdot},M)\}_{n \geq 0}$, we have isomorphisms:
\[
\opn{H}^{r-n+i}(\varDelta^{\centerdot}_{n}(P^{\centerdot},M))
 \cong \opn{H}^{r-n+i}(\varDelta^{\centerdot}_{n+j}(P^{\centerdot},M))
\]
for all $i > 0$ and $\infty \geq j \geq 0$.
\end{lem}

\begin{lem} \label{exttilt1}
For $\{\varDelta^{\centerdot}_{n}(P^{\centerdot},M)\}_{n \geq 0}$ and $\infty \geq n \geq r$, we have
\[
\Hom_{\cat{D}(A)}(P^{\centerdot}, \varDelta^{\centerdot}_{n}(P^{\centerdot}, M)[i]) =0
\]
for all $i \not= r-n$.
\end{lem}

\begin{proof}
Applying $\Hom_{\cat{D}(A)}(P^{\centerdot}, -)$ to $\zeta_{n}(P^{\centerdot}, M)$ ($n \geq 1$),
in case of $0 \leq n \leq r$ we have
\[
\Hom_{\cat{D}(A)}(P^{\centerdot}, \varDelta^{\centerdot}_{n}(P^{\centerdot}, M)[i])=0
\]
for $i > r - n$ or $i < 0$.
Then in case of $n \geq r$ we have
\[
\Hom_{\cat{D}(A)}(P^{\centerdot}, \varDelta^{\centerdot}_{n}(P^{\centerdot}, M)[i])=0
\]
for $i \not= r-n$.
\end{proof}

\begin{thm} \label{exttilt2}
Let $A$ be a symmetric algebra over a field $k$,
and $P^{\centerdot} \in \cat{K}^{\mrm{b}}(\cat{proj}A)$
a partial tilting complex of length $r+1$.
Then the following are equivalent.
\begin{enumerate}
\item $\opn{H}^{i}(\varDelta^{\centerdot}_{r}(P^{\centerdot}, A)) = O$ for all $i > 0$.
\item $\Theta_{n}^{\centerdot}(P^{\centerdot}, A)$ is a tilting complex for any $n \geq r$.
\end{enumerate}
\end{thm}

\begin{proof}
According to the construction of $\varDelta^{\centerdot}_{n}(P^{\centerdot}, A)$,  it is clear that
$\Theta_{n}^{\centerdot}(P^{\centerdot}, A)$ generates $\cat{K}^{\mrm{b}}(\cat{proj}A)$.
By Lemmas \ref{symm} and \ref{exttilt1}, it is easy to see that
$\Theta_{n}^{\centerdot}(P^{\centerdot}, A)$ is a tilting complex for $A$
if and only if $\Hom_{\cat{D}(A)}(\varDelta^{\centerdot}_{n}(P^{\centerdot}, A),
\varDelta^{\centerdot}_{n}(P^{\centerdot}, A)[i])=0$ for all $i > 0$.
By Proposition \ref{exttilt0}, we have 
\[\begin{aligned}
\opn{H}^{i}(\varDelta^{\centerdot}_{r}(P^{\centerdot}, A)) & \cong
\opn{H}^{i}(\varDelta^{\centerdot}_{n}(P^{\centerdot}, A)) \\
& \cong \Hom_{\cat{D}(A)}(A, \varDelta^{\centerdot}_{n}(P^{\centerdot}, A)[i])
\end{aligned}\]
for all $i > 0$.
For $j \leq n$, applying $\Hom_{\cat{D}(A)}(-, \varDelta^{\centerdot}_{n}(P^{\centerdot}, A))$ to 
$\zeta_{j}(P^{\centerdot}, A)$,
we have 
\[
\Hom_{\cat{D}(A)}(\varDelta^{\centerdot}_{j}(P^{\centerdot}, A), 
\varDelta^{\centerdot}_{n}(P^{\centerdot}, A)[i]) \cong 
\Hom_{\cat{D}(A)}(\varDelta^{\centerdot}_{j-1}(P^{\centerdot}, A), 
\varDelta^{\centerdot}_{n}(P^{\centerdot}, A)[i])
\]
for all $i > 0$, because
$\Hom_{\cat{D}(A)}(P^{\centerdot}[j-r-1], 
\varDelta^{\centerdot}_{n}(P^{\centerdot}, A)[i]) = 0$ for all $i \geq 0$.
Therefore $\Hom_{\cat{D}(A)}(A, \varDelta^{\centerdot}_{n}(P^{\centerdot}, A)[i]) = 0$ for all $i > 0$
if and only if $\Hom_{\cat{D}(A)}(\varDelta^{\centerdot}_{n}(P^{\centerdot}, A), \\
\varDelta^{\centerdot}_{n}(P^{\centerdot}, A)[i]) = 0$ for all $i > 0$.
\end{proof}

\begin{cor} \label{exttilt3}
Let $A$ be a symmetric algebra over a field $k$,
$P^{\centerdot} \in \cat{K}^{\mrm{b}}(\cat{proj}A)$
a partial tilting complex of length $r+1$, and $V^{\centerdot}$ the associated bimodule
complex of $P^{\centerdot}$.
Then the following are equivalent.
\begin{enumerate}
\item $\opn{H}^{i}(\varDelta^{\centerdot}_{A}(V^{\centerdot})) = O$ for all $i > 0$.
\item $\Theta_{n}^{\centerdot}(P^{\centerdot}, A)$ is a tilting complex for any $n \geq r$.
\end{enumerate}
\end{cor}

\begin{proof}
According to Corollary \ref{vardelta2}, we have $\varDelta_{A}^{\centerdot}(V^{\centerdot}) \cong
\varDelta_{\infty}^{\centerdot}(P^{\centerdot},A)$ in $\cat{D}(A)$.
Since $\opn{H}^{i}(\varDelta_{\infty}^{\centerdot}(P^{\centerdot},A)) \cong 
\opn{H}^{i}(\varDelta_{r}^{\centerdot}(P^{\centerdot},A))$ for $i > 0$,
we  complete the proof by Theorem \ref{exttilt2}.
\end{proof}

In the case of symmetric algebras, we have a complex version of extensions of classical partial
tilting modules which was showed by Bongartz \cite{Bo}.

\begin{cor} \label{exttilt4}
Let $A$ be a symmetric algebra over a field $k$,
and $P^{\centerdot} \in \cat{K}^{\mrm{b}}(\cat{proj}A)$
a partial tilting complex of length 2.
Then $\Theta_{n}^{\centerdot}(P^{\centerdot}, A)$ is a tilting complex for any $n \geq 1$.
\end{cor}

\begin{proof}
By the construction, $\varDelta_{1}^{i}(P^{\centerdot}, A)=O$
for $i >0$.  According to Theorem \ref{exttilt2} we complete the proof.
\end{proof}

For an object $M$ in an additive category, we denote by $n(M)$ the number of
indecomposable types in $\cat{add}M$.

\begin{cor} \label{exttilt5}
Let $A$ be a symmetric algebra over a field $k$,
and $P^{\centerdot} \in \cat{K}^{\mrm{b}}(\cat{proj}A)$
a partial tilting complex of length 2.
Then the following are equivalent.
\begin{enumerate}
\item $P^{\centerdot}$ is a tilting complex for $A$.
\item $n(P^{\centerdot})=n(A)$.
\end{enumerate}
\end{cor}

\begin{proof}
We may assume $P^{\centerdot}:P^{-1} \to P^{0}$.
Since $\Theta_{1}^{\centerdot}(P^{\centerdot}, A)= 
P^{\centerdot}\oplus\varDelta^{\centerdot}_{1}(P^{\centerdot}, A)$, 
by Corollary \ref{exttilt4}, we have $n(A)=n(\Theta_{1}^{\centerdot}(P^{\centerdot}, A)) =
n(P^{\centerdot})+m$ for some $m \geq 0$.
It is easy to see that $m=0$ if and only if 
$\cat{add}\Theta_{1}^{\centerdot}(P^{\centerdot}, A)= \cat{add}P^{\centerdot}$.
\end{proof}

\begin{lem} \label{ddotadj2}
Let $\theta:\bsym{1}_{\cat{D}(eAe)} \to j^{e *}_{A}j^{e}_{A !}$ be the adjunction arrow,
and let $X^{\centerdot} \in \cat{D}(eAe)$ and $Y^{\centerdot} \in \cat{D}(A)$.
For $h \in \Hom_{\cat{D}(A)}(j^{e}_{A !}(X^{\centerdot}),Y^{\centerdot})$,
let $\Phi (h)= j^{e *}_{A}(h) \circ \theta_X$, then
$\Phi : \Hom_{\cat{D}(A)}(j^{e}_{A !}(X^{\centerdot}),Y^{\centerdot}) \xarr{\sim}
\Hom_{\cat{D}(A)}(X^{\centerdot},j^{e *}_{A}Y^{\centerdot})$
is an isomorphism as $\opn{End}_{\cat{D}(A)}(X^{\centerdot})$-modules.
\end{lem}

\begin{thm} \label{extrecoll}
Let $A$ be a symmetric algebra over a field $k$, $e$ an idempotent of $A$, 
$Q^{\centerdot} \in \cat{K}^{\mrm{b}}(\cat{proj}eAe)$ a tilting complex for $eAe$,
and $P^{\centerdot}=j^{e}_{A !}(Q^{\centerdot}) \in \cat{K}^{\mrm{b}}(\cat{proj}A)$ with
$l(P^{\centerdot})=r+1$.
For $n \geq r$, the following hold.
\begin{enumerate}
\item $\Theta_{n}^{\centerdot}(P^{\centerdot},A)$ is
a recollement tilting complex related to $e$.
\item $A/AeA \cong B/BfB$, where 
$B=\opn{End}_{\cat{D}(A)}(\Theta_{n}^{\centerdot}(P^{\centerdot},A))$ and $f$ is an
idempotent of $B$ corresponding to $e$.
\end{enumerate}
\end{thm}

\begin{proof}
We may assume $P^{\centerdot}:P^{-r} \to \ldots P^{-1} \to P^{0}$.
Since $j^{e}_{A !}$ is fully faithful, 
$\Hom_{\cat{D}(A)}(P^{\centerdot}, P^{\centerdot}[i])=0$ for $i \not= 0$.  
Consider a family $\{\varDelta^{\centerdot}_{n}(P^{\centerdot},A)\}_{n\geq 0}$ 
of Definition \ref{consttilt} and triangles $\zeta_{n}(P^{\centerdot}, A)$:
\[
P^{\centerdot}_{n}[n-r-1] \xarr{g_n} \varDelta^{\centerdot}_{n-1}(P^{\centerdot},A)
\xarr{h_n} \varDelta^{\centerdot}_{n}(P^{\centerdot},A) \arr P^{\centerdot}_{n}[n-r] .
\]
The morphism $\Phi$ of Lemma \ref{ddotadj2} induces isomorphisms between exact sequences in $\cat{Mod}B$:
{\small \[\begin{array}{cc}
\Hom_{\cat{D}(A)}(P^{\centerdot}, P^{\centerdot}_{n}[n-r-1+i]) \to 
& \Hom_{\cat{D}(A)}(P^{\centerdot},\varDelta^{\centerdot}_{n-1}(P^{\centerdot},A)[i]) \to \\
\downarrow \Phi & \downarrow \Phi \\
\Hom_{\cat{D}(eAe)}(Q^{\centerdot}, j^{e *}_{A}P^{\centerdot}_{n}[n-r-1+i]) \to 
& \Hom_{\cat{D}(eAe)}(Q^{\centerdot},j^{e *}_{A}\varDelta^{\centerdot}_{n-1}(P^{\centerdot},A)[i]) \to \\
\\
\Hom_{\cat{D}(A)}(P^{\centerdot},\varDelta^{\centerdot}_{n}(P^{\centerdot},A)[i]) \to 
& \Hom_{\cat{D}(A)}(P^{\centerdot}, P^{\centerdot}_{n}[n-r+i])\\
\downarrow \Phi &\downarrow \Phi \\
\Hom_{\cat{D}(eAe)}(Q^{\centerdot},j^{e *}_{A}\varDelta^{\centerdot}_{n}(P^{\centerdot},A)[i]) \to
& \Hom_{\cat{D}(eAe)}(Q^{\centerdot}, j^{e *}_{A}P^{\centerdot}_{n}[n-r+i])
\end{array}\]}\par\noindent
for all $i$.  By Lemma \ref{ddotadj2}, we have $j^{e *}_{A}(\zeta_{n}(P^{\centerdot}, A)) \cong 
\zeta_{n}(Q^{\centerdot}, j^{e *}_{A}A)$ in $\cat{D}(eAe)$, and then
$\{j^{e *}_{A}(\varDelta^{\centerdot}_{n}(P^{\centerdot},A))\}_{n \geq 0}
\cong \{\varDelta^{\centerdot}_{n}(Q^{\centerdot},Ae)\}_{n \geq 0}$.
By lemma \ref{exttilt1}, it is easy to see that
\[
\Hom_{\cat{D}(eAe)}(Q^{\centerdot}, \varDelta^{\centerdot}_{\infty}(Q^{\centerdot},Ae)[i]) =0
\]
for all $i \in \mbb{Z}$.
Since $Q^{\centerdot}$ is a tilting complex for $eAe$,
$\varDelta^{\centerdot}_{\infty}(Q^{\centerdot},Ae)$ is a null complex, that is 
$\opn{H}^{i}(\varDelta^{\centerdot}_{\infty}(Q^{\centerdot},Ae)) =O$ for all $i \in \mbb{Z}$.
By Lemma \ref{exttilt0}, for $n \geq r$ we have
$\opn{H}^{i}(\varDelta^{\centerdot}_{n}(Q^{\centerdot},Ae))=O$ for all $i > 0$.
By the above isomorphism, for $n \geq r$ we have
$\opn{H}^{i}(\varDelta^{\centerdot}_{n}(P^{\centerdot},A)) \in \cat{Mod}A/AeA$ for all $i > 0$.
On the other hand, $\varDelta^{\centerdot}_{n}(P^{\centerdot},A)$ has the form:
\[
R^{\centerdot}:R^{-n} \to \ldots \to R^{0} \to R^{1} \to \ldots \to R^{r-1} ,
\]
where $R^{i} \in \cat{add}eA$ for $i \not= 0$, and $R^{0}=A\oplus R^{' 0}$
with $R^{' 0} \in \cat{add}eA$.
Since $\Hom_{A}(eA, \cat{Mod}A/AeA)=0$, it is easy to see that
$\varDelta^{\centerdot}_{n}(P^{\centerdot},A) \cong 
\sigma_{\leq 0}\varDelta^{\centerdot}_{n}(P^{\centerdot},A)$
($\cong \sigma_{\leq 0}\ldots\sigma_{\leq r-2}\varDelta^{\centerdot}_{n}(P^{\centerdot},A)$
if $r \geq 2$).
Therefore, $\opn{H}^{i}(\varDelta^{\centerdot}_{n}(P^{\centerdot},A))=O$ for
all $i > 0$, and hence $\Theta_{n}^{\centerdot}(P^{\centerdot},A)$ is
a recollement tilting complex related to $e$ by Theorem \ref{exttilt2}.
Since $\Theta_{n}^{\centerdot}(P^{\centerdot},A) \cong P^{\centerdot}[n-r]\oplus R^{\centerdot}$
and $j^{e}_{A !}(X^{\centerdot})\lten_{A}\varDelta^{\centerdot}_{A}(e)=
i^{e *}_{A}j^{e}_{A !}(X^{\centerdot})=O$ for $X^{\centerdot} 
\in \cat{D}(eAe)$, we have an isomorphism
$\Theta_{n}^{\centerdot}(P^{\centerdot},A)\dten_{A}\varDelta^{\centerdot}_{A}(e) \cong
\varDelta^{\centerdot}_{A}(e)$ in $\cat{D}(A)$.
By Proposition \ref{Moritaeqv1}, we complete the proof.
\end{proof}

\begin{cor} \label{Moritaeqv2}
Under the condition Theorem \ref{extrecoll}, let ${}_{B}T^{\centerdot}_{A}$ be 
the associated two-sided tilting complex of $\Theta_{n}^{\centerdot}(P^{\centerdot},A)$.
Then the standard equivalence $\rhom_{A}(T^{\centerdot}, -):\cat{D}(A) \xarr{\sim} \cat{D}(B)$ induces
an equivalence $R^{0}\dhom_{A}(T^{\centerdot}, -)|_{\cat{Mod}A/AeA}:
\cat{Mod}A/AeA \xarr{\sim} \cat{Mod}B/BfB$.
\end{cor}

\begin{proof}
By the proof of Theorem \ref{extrecoll}, we have 
$T^{\centerdot}\lten_{A}\varDelta_{A}^{\centerdot}(e) \cong \varDelta_{A}^{\centerdot}(e)$
in $\cat{D}(A)$.
By Proposition \ref{Moritaeqv1}, we complete the proof.
\end{proof}

\begin{rem} \label{rem1}
For a symmetric algebra $A$ over a field $k$ and an idempotent $e$ of $A$, $eAe$ is also a symmetric
$k$-algebra.  Therefore, we have constructions of tilting complexes with respect to any sequence
of idempotents of $A$.
Moreover, if a recollement $\{\cat{D}_{A/AeA}(A), \cat{D}(A), \cat{D}(eAe)\}$
is triangle equivalent to a recollement $\{\cat{D}_{B/BfB}(B), \\ \cat{D}(B), \cat{D}(fBf)\}$, then
$B$ and $fBf$ are also symmetric $k$-algebras.
\end{rem}

\begin{rem} \label{rem2}
According to \cite{Rd2}, under the condition of Theorem \ref{extrecoll}
we have a stable equivalence $\underline{\cat{mod}}A \xarr{\sim} \underline{\cat{mod}}B$
which sends $A/AeA$-modules to $B/BfB$-modules, where
$\underline{\cat{mod}}A, \underline{\cat{mod}}B$ are stable categories
of finitely generated modules.
In particular, this equivalence sends simple $A/AeA$-modules to simple $B/BfB$-modules.
\end{rem}

\begin{rem} \label{rem3}
Let $A$ be a ring, and $e$ an idempotent of $A$ such that
there is a finitely generated projective resolution of $Ae$ in $\cat{Mod}eAe$.
Then Hoshino and Kato showed that $\Theta_{n}^{\centerdot}(eA,A)$
is a tilting complex if and only if $\opn{Ext}_{A}^{i}(A/AeA,eA)=0$ for $0 \leq i < n$
(\cite{HK2}).
In even this case, we have also $A/AeA \cong B/BfB$, where 
$B=\End_{\cat{D}(A)}(\Theta_{n}^{\centerdot}(eA,A))$ and $f$ is an idempotent of $B$
corresponding to $e$.
Moreover if $A$, $B$ are projective algebras over a commutative ring $k$, then
by Proposition \ref{Moritaeqv1}
the standard equivalence induces an equivalence $\cat{Mod}A/AeA \xarr{\sim} \cat{Mod}B/BfB$.
\end{rem}


\end{document}